\documentclass{amsart}
\usepackage{mathptmx}
\usepackage{pxfonts}
\usepackage{amssymb}
\usepackage{amsfonts}
\usepackage{amsmath}
\usepackage{graphicx}
\usepackage{shadow}
\usepackage{color}
\usepackage{enumitem}
\usepackage[all]{xy}
\newtheorem{theorem}{Theorem}[section]
\newtheorem{corollary}[theorem]{Corollary}
\newtheorem{lemma}[theorem]{Lemma}
\newtheorem{proposition}[theorem]{Proposition}

\theoremstyle{definition}
\newtheorem{definition}[theorem]{Definition}
\newtheorem{example}[theorem]{Example}

\theoremstyle{remark}
\newtheorem{remark}[theorem]{Remark}

\newcommand{\field}[1]{\mathbb{#1}}

\newcommand{\R}{\field{R}}

\DeclareMathOperator{\htt}{ht}

\DeclareMathOperator{\td}{t}
\DeclareMathOperator{\t.d.}{t.d.}

\DeclareMathOperator{\Ker}{Ker}

\DeclareMathOperator{\Spec}{Spec}
\DeclareMathOperator{\Max}{Max}

\begin{document}
%\baselineskip=14pt
%%%%%%%%%%%%%%%%%%%%%%%%%%%%%%%%%%%%%%%%%%%%%%%%%%%%%%%%%%%%%%%%%%%%%%%%%%%%%%%%%%%%%%%%%%%%%%%%%%%%%%%%%%%%%%%%%%%%%%%%%%%%%%%%%%%%%%%%%%%%%%
%%%%%%%%%%%%%%%%%%%%%%%%%%%%%%%%%%%%%%%%%%%%%%%%%%%%%%%%%%%%%%%%%%%%%%%%%%%%%%%%%%%%%%%%%%%%%%%%%%%%%%%%%%%%%%%%%%%%%%%%%%%%%%%%%%%%%%%%%%%%%%
%%%%%%%%%%%%%%%%%%%%%%%%%%%%%%%%%%%%%%%%%%%%%%%%%%%%%%%%%%%%%%%%%%%%%%%%%%%%%%%%%%%%%%%%%%%%%%%%%%%%%%%%%%%%%%%%%%%%%%%%%%%%%%%%%%%%%%%%%%%%%%
%%%%%%%%%%%%%%%%%%%%%%%%%%%%%%%%%%%%%%%%%%%%%%%%%%%%%%%%%%%%%%%%%%%%%%%%%%%%%%%%%%%%%%%%%%%%%%%%%%%%%%%%%%%%%%%%%%%%%%%%%%%%%%%%%%%%%%%%%%%%%%

\title[Dimension theory of tensor products of AF-rings]{Dimension theory of tensor products of AF-rings}

\author[Salah Kabbaj]{Salah Kabbaj}
\address{Department of Mathematics and Statistics, King Fahd University of Petroleum and Minerals, Dhahran 31261, Saudi Arabia}
\email{kabbaj@kfupm.edu.sa}

\date{\today}

\subjclass[2010]{13C15, 13B24, 13F05, 13H05, 13F20, 13B30, 13E05, 13D05}

\keywords{Tensor product of algebras, Krull dimension, valuative dimension, Jaffard ring, locally Jaffard ring, altitude formula, AF-ring, pullback}

%\dedicatory{}

%%%%%%%%%%%%%%%%%%%%%%%%%%%%%%%%%%%%%%%%%%%%%%%%%%%%%%%%%%%%%%%%%%%%%%%%%%%%%%%%%%%%%%%%%%%%%%%%%%%%%%%%%%%%%%%%%%%%%%%%%%%%%%%%%%%%%%%%%%%%%%
%%%%%%%%%%%%%%%%%%%%%%%%%%%%%%%%%%%%%%%%%%%%%%%%%%%%%%%%%%%%%%%%%%%%%%%%%%%%%%%%%%%%%%%%%%%%%%%%%%%%%%%%%%%%%%%%%%%%%%%%%%%%%%%%%%%%%%%%%%%%%%
%%%%%%%%%%%%%%%%%%%%%%%%%%%%%%%%%%%%%%%%%%%%%%%%%%%%%%%%%%%%%%%%%%%%%%%%%%%%%%%%%%%%%%%%%%%%%%%%%%%%%%%%%%%%%%%%%%%%%%%%%%%%%%%%%%%%%%%%%%%%%%
%%%%%%%%%%%%%%%%%%%%%%%%%%%%%%%%%%%%%%%%%%%%%%%%%%%%%%%%%%%%%%%%%%%%%%%%%%%%%%%%%%%%%%%%%%%%%%%%%%%%%%%%%%%%%%%%%%%%%%%%%%%%%%%%%%%%%%%%%%%%%%
\begin{abstract}
AF-rings are algebras over a field $k$ which satisfy the Altitude Formula over $k$. This paper surveys a few works in the literature on the Krull and valuative dimensions of tensor products of AF-rings. The first section extends Wadsworth's classical results on the Krull dimension of AF-domains to the larger class of AF-rings. It also provides formulas for computing the valuative dimension with effect on the transfer of the (locally) Jaffard property. The second section studies tensor products of AF-rings over a zero-dimensional ring. Most results on algebras over a field are extended to these general constructions.  The third section establishes formulas for the Krull and valuative dimensions of tensor products of pullbacks issued from AF-domains. Throughout, examples are provided to illustrate the scope and limits of the results.
\end{abstract}
\maketitle

%%%%%%%%%%%%%%%%%%%%%%%%%%%%%%%%%%%%%%%%%%%%%%%%%%%%%%%%%%%%%%%%%%%%%%%%%%%%%%%%%%%%%%%%%%%%%%%%%%%%%%%%%%%%%%%%%%%%%%%%%%%%%%%%%%%%%%%%%%%%%%
%%%%%%%%%%%%%%%%%%%%%%%%%%%%%%%%%%%%%%%%%%%%%%%%%%%%%%%%%%%%%%%%%%%%%%%%%%%%%%%%%%%%%%%%%%%%%%%%%%%%%%%%%%%%%%%%%%%%%%%%%%%%%%%%%%%%%%%%%%%%%%
%%%%%%%%%%%%%%%%%%%%%%%%%%%%%%%%%%%%%%%%%%%%%%%%%%%%%%%%%%%%%%%%%%%%%%%%%%%%%%%%%%%%%%%%%%%%%%%%%%%%%%%%%%%%%%%%%%%%%%%%%%%%%%%%%%%%%%%%%%%%%%
%%%%%%%%%%%%%%%%%%%%%%%%%%%%%%%%%%%%%%%%%%%%%%%%%%%%%%%%%%%%%%%%%%%%%%%%%%%%%%%%%%%%%%%%%%%%%%%%%%%%%%%%%%%%%%%%%%%%%%%%%%%%%%%%%%%%%%%%%%%%%%
\begin{center}\begin{minipage}{0.75\textwidth}
{\bf\footnotesize
\tableofcontents
}
\end{minipage}\end{center}
\bigskip

%\newpage
%%%%%%%%%%%%%%%%%%%%%%%%%%%%%%%%%%%%%%%%%%%%%%%%%%%%%%%%%%%%%%%%%%%%%%%%%%%%%%%%%%%%%%%%%%%%%%%%%%%%%%%%%%%%%%%%%%%%%%%%%%%%%%%%%%%%%%%%%%%%%%
%%%%%%%%%%%%%%%%%%%%%%%%%%%%%%%%%%%%%%%%%%%%%%%%%%%%%%%%%%%%%%%%%%%%%%%%%%%%%%%%%%%%%%%%%%%%%%%%%%%%%%%%%%%%%%%%%%%%%%%%%%%%%%%%%%%%%%%%%%%%%%
%%%%%%%%%%%%%%%%%%%%%%%%%%%%%%%%%%%%%%%%%%%%%%%%%%%%%%%%%%%%%%%%%%%%%%%%%%%%%%%%%%%%%%%%%%%%%%%%%%%%%%%%%%%%%%%%%%%%%%%%%%%%%%%%%%%%%%%%%%%%%%
%%%%%%%%%%%%%%%%%%%%%%%%%%%%%%%%%%%%%%%%%%%%%%%%%%%%%%%%%%%%%%%%%%%%%%%%%%%%%%%%%%%%%%%%%%%%%%%%%%%%%%%%%%%%%%%%%%%%%%%%%%%%%%%%%%%%%%%%%%%%%%
\section*{Introduction}\label{i}

\noindent Throughout this paper, all rings and, for a given field $k$, all $k$-algebras are assumed to be commutative with identity element and \emph{have finite transcendence degree over $k$}. For a ring $A$, we shall use  $\Spec(A)$ and $\Max(A)$ to denote, respectively, the sets of all prime ideals and maximal ideals of $A$. Also, we will denote by $A[n]$ the polynomial ring $A[X_1,\dots,X_n]$ and by $p[n]$ the prime ideal $p[X_1,\dots,X_n]$ in $A[n]$, for any $p\in\Spec(A)$ and positive integer $n$.

A finite-dimensional domain $R$ is said to be Jaffard if $\dim(R[n])= n + \dim(R)$ for all $n \geq 1$; equivalently, if $\dim(R) = \dim_{v}(R)$, where $\dim(R)$ denotes the Krull dimension of $R$ and $\dim_{v}(R)$ denotes its valuative dimension (i.e., the supremum of dimensions of the valuation overrings of $R$). Since this notion does not carry over to localizations, $R$ is said to be locally Jaffard if $R_p$ is Jaffard for each $p\in\Spec(R)$ (equivalently,  $\htt(p[n])=\htt(p)$, $\forall\ p\in\Spec(A)$). The class of Jaffard domains contains most of the well-known classes of rings involved in dimension theory such as Noetherian domains, Pr\"ufer domains, universally catenarian domains, and stably strong S-domains. Analogous definitions are given in Cahen's paper \cite{Cahen} for a finite-dimensional arbitrary ring (i.e., possibly with zero-divisors). We assume familiarity with these concepts as in \cite{ABDFK,ADKM,BDK,BDF,BK,BMRH,DFK,FK1990,FK,J,K1,K2,MM}, and any unreferenced material is standard as in \cite{G1,Ka,Ma}.

Let $k$ be a field and let $\td(A)$ denote the transcendence degree over $k$ of a $k$-algebra $A$. If $A$ is not a domain, then by definition $\td(A):=\max\big\{\td(A/p)\mid p\in\Spec(A)\big\}$.

%%%%%%%%%%%%%%%%%%%%%%%%%%%%%%%%%%%%%%%%%%%%%%%%%%%%%%%%%%%%%%%%%
%%%%%%%%%%%%%%%%%%%%%%%%%%%%%%%%%%%%%%%%%%%%%%%%%%%%%%%%%%%%%%%%%
\begin{definition}
A $k$-algebra $A$ is an AF-ring if it satisfies the Altitude Formula over $k$; that is, $\htt(p) +\td(A/p)=\td(A_p)$, for each $p\in\Spec(A)$.
\end{definition}

Examples and basic properties of AF-rings are provided at the beginning of the next section. In 1977, Sharp proved in \cite{Sh} that
$$\dim(K_1\otimes_kK_2)=\min\big(\td(K_1),\td(K_2)\big)$$
for any field extensions $K_1$ and $K_2$ of $k$. In 1978, Sharp and Vamos generalized this result to the tensor product of a finite number of field extensions of $k$ \cite{ShVa}. In 1979, Wadsworth extended their results on field extensions to the larger class of AF-domains \cite{W}; namely, he proved that if $D_1$ and $D_2$ are AF-domains, then
$$\dim(D_1\otimes_kD_2) =\min\big(\dim(D_1) + \td(D_2), \td(D_1) + \dim(D_2)\big).$$
Moreover, he established a formula for $\dim(D\otimes_kR)$ which holds for an AF-domain $D$, with no  restriction on the ring $R$. He also proved that for any prime ideal $p$ of an AF-ring $A$ and, for any $n\ge1$, $\htt(p[n])=\htt(p)$ (i.e., $A$ is locally Jaffard). In \cite{Giro}, Girolami studied the class of AF-domains with respect to the class of $k$-algebras which are stably strong $S$-domains and examined the behavior of the notion of AF-domain for certain pullback constructions. An upper bound was then given for the valuative  dimension of the tensor product of two $k$-algebras; more exactly, if $A_1$ and $A_2$ are $k$-algebras, then $$\dim_v(A_1\otimes_kA_2)\le\min\big(\dim_v(A_1) +\td(A_2),\td(A_1) + \dim_v(A_2)\big).$$

This paper surveys a few works in the literature on the Krull and valuative dimensions of tensor products of AF-rings. The first section extends Wadsworth's classical results on the Krull dimension of AF-domains to the larger class of AF-rings. It also provides formulas for computing the valuative dimension with effect on the transfer of the (locally) Jaffard property. The second section studies tensor products of AF-rings over a zero-dimensional ring. Most results on algebras over a field are extended to these general constructions.  The third section establishes formulas for the Krull and valuative dimensions of tensor products of pullbacks issued from AF-domains. Throughout, examples are provided to illustrate the scope and limits of the results.

The three main papers involved in this survey are \cite{BGK1,BGK2,BGK3}, which were co-authored with Samir Bouchiba (University of Meknes) and Florida Girolami (University of Rome) and published in 1997/1999.

%\newpage
%%%%%%%%%%%%%%%%%%%%%%%%%%%%%%%%%%%%%%%%%%%%%%%%%%%%%%%%%%%%%%%%%%%%%%%%%%%%%%%%%%%%%%%%%%%%%%%%%%%%%%%%%%%%%%%%%%%%%%%%%%%%%%%%%%%%%%%%%%%%%%
%%%%%%%%%%%%%%%%%%%%%%%%%%%%%%%%%%%%%%%%%%%%%%%%%%%%%%%%%%%%%%%%%%%%%%%%%%%%%%%%%%%%%%%%%%%%%%%%%%%%%%%%%%%%%%%%%%%%%%%%%%%%%%%%%%%%%%%%%%%%%%
%%%%%%%%%%%%%%%%%%%%%%%%%%%%%%%%%%%%%%%%%%%%%%%%%%%%%%%%%%%%%%%%%%%%%%%%%%%%%%%%%%%%%%%%%%%%%%%%%%%%%%%%%%%%%%%%%%%%%%%%%%%%%%%%%%%%%%%%%%%%%%
%%%%%%%%%%%%%%%%%%%%%%%%%%%%%%%%%%%%%%%%%%%%%%%%%%%%%%%%%%%%%%%%%%%%%%%%%%%%%%%%%%%%%%%%%%%%%%%%%%%%%%%%%%%%%%%%%%%%%%%%%%%%%%%%%%%%%%%%%%%%%%
\section{Tensor products of AF-rings over a field}\label{af}

\noindent  This section is devoted to \cite{BGK1}. First it extends some classical results (on the Krull dimension) known for the class of AF-domains to the class of AF-rings over a field. Then it provides formulas for computing the valuative dimension of tensor products emanating from AF-rings with effect on the possible transfer of the notion of (locally) Jaffard ring to  these constructions.

Throughout this section $k$ will denote a field and by a ring we mean a $k$-algebra. Also, algebras (resp., tensor products) are taken over (resp., relative to) $k$. For the reader's convenience, we first recall some basic properties of AF-rings.

%%%%%%%%%%%%%%%%%%%%%%%%%%%%%%%%%%%%%%%%%%%%%%%%%%%%%%%%%%%%%%%%%
%%%%%%%%%%%%%%%%%%%%%%%%%%%%%%%%%%%%%%%%%%%%%%%%%%%%%%%%%%%%%%%%%
\begin{remark}[\cite{Giro,W}]
Let $\mathcal{A}$ denote the class of AF-rings over $k$ and let $n$ be a positive integer. Then:
\begin{enumerate}
\item Any finitely generated algebra and its integral extensions belong to $\mathcal{A}$.
\item If $A\in\mathcal{A}$, then $S^{-1}A\in\mathcal{A}$, for every multiplicative subset $S$ of $A$.
\item If $A\in\mathcal{A}$, then $A[n]\in\mathcal{A}$.
\item If  $A_1,\dots, A_n\in\mathcal{A}$, then $A_1\otimes\dots\otimes A_n\in\mathcal{A}$.
\item If  $A_1,\dots, A_n\in\mathcal{A}$, then $A_1\times\dots\times A_n\in\mathcal{A}$.
\item If $A\in\mathcal{A}$, then $A$ is locally Jaffard.
\item The class $\mathcal{A}$ is not stable under factor rings. However, if $A$ is a catenarian AF-domain, then $A/p\in\mathcal{A},\ \forall\ p\in\Spec(A)$.
\end{enumerate}
\end{remark}

%%%%%%%%%%%%%%%%%%%%%%%%%%%%%%%%%%%%%%%%%%%%%%%%%%%%%%%%%%%%%%%%%%%%%%%%%%%%%%%%%%%%%%%%%%%%%%%%%%%%%
%%%%%%%%%%%%%%%%%%%%%%%%%%%%%%%%%%%%%%%%%%%%%%%%%%%%%%%%%%%%%%%%%%%%%%%%%%%%%%%%%%%%%%%%%%%%%%%%%%%%%
%%%%%%%%%%%%%%%%%%%%%%%%%%%%%%%%%%%%%%%%%%%%%%%%%%%%%%%%%%%%%%%%%%%%%%%%%%%%%%%%%%%%%%%%%%%%%%%%%%%%%
%%%%%%%%%%%%%%%%%%%%%%%%%%%%%%%%%%%%%%%%%%%%%%%%%%%%%%%%%%%%%%%%%%%%%%%%%%%%%%%%%%%%%%%%%%%%%%%%%%%%%
\subsection{Krull dimension}

\noindent This subsection aims at extending Wadsworth's results on AF-domains to the class of AF-rings. The first technical result links the transcendence degree of a localization of a tensor product to the transcendence degrees of its respective components.

%%%%%%%%%%%%%%%%%%%%%%%%%%%%%%%%%%%%%%%%%%%%%%%%%%%%%%%%%%%%%%%%%
%%%%%%%%%%%%%%%%%%%%%%%%%%%%%%%%%%%%%%%%%%%%%%%%%%%%%%%%%%%%%%%%%
\begin{lemma}\label{af-lem1}
 Let  $A_1,\dots, A_n$ be AF-rings and let $P\in \Spec(A_1\otimes\dots\otimes A_n)$. Then
$$\td\big((A_1\otimes\dots\otimes A_n)_{P}\big) = \sum_{1\leq i\leq n}\td(A_{i_{p_{i}}})$$
where $p_i := P\cap A_i$ for $i=1, \dots, n$.
\end{lemma}

As an immediate consequence of this lemma, we obtain the following known result for AF-domains.

%%%%%%%%%%%%%%%%%%%%%%%%%%%%%%%%%%%%%%%%%%%%%%%%%%%%%%%%%%%%%%%%%
%%%%%%%%%%%%%%%%%%%%%%%%%%%%%%%%%%%%%%%%%%%%%%%%%%%%%%%%%%%%%%%%%
\begin{corollary}
 Let  $D_1,\dots, D_n$ be AF-domains and let $P\in\Spec(D_1\otimes\dots\otimes D_n)$. Then
$$\td\big((D_1\otimes\dots\otimes D_n)_{P}\big) = \td(D_1\otimes\dots\otimes D_n)= \sum_{1\leq i\leq n}\td(D_{i}).$$
\end{corollary}

The following simple statement has important consequences on some of the following results.

%%%%%%%%%%%%%%%%%%%%%%%%%%%%%%%%%%%%%%%%%%%%%%%%%%%%%%%%%%%%%%%%%
%%%%%%%%%%%%%%%%%%%%%%%%%%%%%%%%%%%%%%%%%%%%%%%%%%%%%%%%%%%%%%%%%
\begin{lemma}
Let $A$ be an AF-ring and $p\in\Spec(A)$. Let $p_{o}$ be a minimal prime ideal of $A$ contained in $p$ such that $\htt(p) = \htt(p/p_{o})$. Then, $\td(A_{p})=\td(A_{p_{o}})$.
\end{lemma}

In order to proceed with the main results, let us recall from \cite{W} the following functions: Given two rings $A$ and $B$ with $p\in \Spec(A)$ and $q\in\Spec(B)$, consider the function
$$\delta(p, q) = \max\big\{\htt(P)\mid P\in\Spec(A\otimes B)\ \text{with}\ P\cap A = p\ \text{and}\ P\cap B = q\big\}.$$
Given a ring $A$, $p\in\Spec(A)$ and $d,s$ integers with $0\leq d\leq s$, consider the two functions
$$\Delta(s, d, p) = \htt(p[s]) + \min\big(s, d + \td(A/p)\big)$$
$$D(s, d, A) =\max\big\{\Delta(s, d, p)\mid p\in\Spec(A)\big\}.$$

The main result of this section provides a formula for the Krull dimension of a tensor product.

%%%%%%%%%%%%%%%%%%%%%%%%%%%%%%%%%%%%%%%%%%%%%%%%%%%%%%%%%%%%%%%%%
%%%%%%%%%%%%%%%%%%%%%%%%%%%%%%%%%%%%%%%%%%%%%%%%%%%%%%%%%%%%%%%%%
\begin{theorem}\label{af-main1}
Let $A$ be an AF-ring and $B$ an arbitrary ring. Then:
\begin{enumerate}
\item $\delta(p, q) = \Delta(\td(A_{p}), \htt(p), q)$, for any $p\in\Spec(A)$ and $q\in\Spec(B)$.
\item $\dim(A\otimes B) = \max\big\{D(\td(A_{p}), \htt(p), B)\mid p\in\Spec(A)\big\}$.
\end{enumerate}
\end{theorem}

Notice that (1) is the most important part of the above theorem. Its proof relies on the above two lemmas after reduction -via localization techniques- to the case where $B$ is a field. Then, the result upon $\dim(A\otimes B)$ derives directly from the definitions of $\delta$, $\Delta$, and $D$.

In case both $A$ and $B$ are AF-rings, we get the following more explicit formula for the Krull dimension.

%%%%%%%%%%%%%%%%%%%%%%%%%%%%%%%%%%%%%%%%%%%%%%%%%%%%%%%%%%%%%%%%%
%%%%%%%%%%%%%%%%%%%%%%%%%%%%%%%%%%%%%%%%%%%%%%%%%%%%%%%%%%%%%%%%%
\begin{corollary}\label{af-cor1}
Let $A$ and $B$ be two AF-rings. Then:
$$\dim(A\otimes B) = \max\big\{\min\big(\htt(p)+\td(A_{q}), \td(A_{p})+\htt(q)\big)\mid p\in\Spec(A), q\in\Spec(B)\big\}.$$
\end{corollary}

The general case of $n$ AF-rings ($n\geq 2$) can be proved by induction on $n$ via Corollary~\ref{af-cor1} and Lemma~\ref{af-lem1}. Namely, we have the following result.

%%%%%%%%%%%%%%%%%%%%%%%%%%%%%%%%%%%%%%%%%%%%%%%%%%%%%%%%%%%%%%%%%
%%%%%%%%%%%%%%%%%%%%%%%%%%%%%%%%%%%%%%%%%%%%%%%%%%%%%%%%%%%%%%%%%
\begin{corollary}
 Let  $A_1,\dots, A_n$ be AF-rings. Then:
$$\dim(A_{1}\otimes\dots \otimes A_{n}) = \max\big\{\min\big(\htt(p_{i})+\sum_{j\not=i}\td(A_{p_{j}})\big)_{1\leq i\leq n}\mid p_{i}\in\Spec(A_{i})\big\}.$$
\end{corollary}

Notice that $D(s,d, A)$  is a nondecreasing function of the first two  arguments and, hence,
one can restrict the formulas in the above three results to the maximal ideals.

Wadsworth's well-known result \cite[Theorem 3.8]{W} on the Krull dimension of the tensor product of $n$ AF-domains reads as follows:  Let  $D_1,\dots, D_n$ be AF-domains with $n\geq 2$. Then
$$\dim(D_{1}\otimes\dots \otimes D_{n}) = \sum_{1\leq i\leq n}\td(D_{i})\ -\ \max\big\{\td(D_{i}) - \dim(D_{i})\big\}_{1\leq i\leq n}.$$

This formula does not hold in general for AF-rings, as shown by the following example.

%%%%%%%%%%%%%%%%%%%%%%%%%%%%%%%%%%%%%%%%%%%%%%%%%%%%%%%%%%%%%%%%%
%%%%%%%%%%%%%%%%%%%%%%%%%%%%%%%%%%%%%%%%%%%%%%%%%%%%%%%%%%%%%%%%%
\begin{example}
Let $R_1:= k[X_1,X_2,X_3]_{(X_1)}$, $R_2 :=  k[X_1,X_2]$, $A_1:= R_1\times R_2$,  and $A_2 := k[X_1,X_2]_{(X_1)}$. Clearly, $A_1$  is an AF-ring with $\dim(A_1) = 2$  and $t(A_1)= 3$; and $A_2$  is an AF-domain with  $\dim(A_2) = 1$  and  $\td(A_2) = 2$.  By Corollary~\ref{af-cor1}, one can check that $\dim(A_1\otimes A_2) = 3\lneqq \td(A_1) + \td(A_2) - 1 = 4$.
\end{example}

The second main result of this section establishes necessary and sufficient conditions for a tensor product of AF-rings to satisfy Wadsworth's aforementioned formula.

%%%%%%%%%%%%%%%%%%%%%%%%%%%%%%%%%%%%%%%%%%%%%%%%%%%%%%%%%%%%%%%%%
%%%%%%%%%%%%%%%%%%%%%%%%%%%%%%%%%%%%%%%%%%%%%%%%%%%%%%%%%%%%%%%%%
\begin{theorem}\label{af-main2}
Let  $A_1,\dots, A_n$ be AF-rings. Then, the following assertions are equivalent:
\begin{enumerate}[label=(\rm\roman*)]
\item $\dim(A_{1}\otimes\dots \otimes A_{n}) = \sum_{1\leq i\leq n}\td(A_{i})\ -\ \max\big\{\td(A_{i}) - \dim(A_{i})\big\}_{1\leq i\leq n}$;
\item For each $i =1,\dots, n$  there exists $M_i\in\Max(A_i)$ such that $\htt(M_{i_{o}})=\dim(A_{i_{o}})$ for some $i_{o}\in\{1,..., n\}$ and, for all $i\not=i_{o}$, $\td(A_{iM_i}) = \td(A_{i})$ \& $\td(A_i/M_i) \leq\td(A_{i_{o}}/M_{i_{o}})$.
\end{enumerate}
\end{theorem}

Next, we give some applications of this result.

%%%%%%%%%%%%%%%%%%%%%%%%%%%%%%%%%%%%%%%%%%%%%%%%%%%%%%%%%%%%%%%%%
%%%%%%%%%%%%%%%%%%%%%%%%%%%%%%%%%%%%%%%%%%%%%%%%%%%%%%%%%%%%%%%%%
\begin{corollary}
 Let  $A_1,\dots, A_n$ be AF-rings such that, for each $i =1,\dots, n$,  there exists $M_i\in\Max(A_i)$ with $\htt(M_{i})=\dim(A_{i})$ and $\td(A_{iM_i}) = \td(A_{i})$. Then  $$\dim(A_{1}\otimes\dots \otimes A_{n}) = \sum_{1\leq i\leq n}\td(A_{i})\ -\ \max\big\{\td(A_{i}) - \dim(A_{i})\big\}_{1\leq i\leq n}.$$
\end{corollary}

%%%%%%%%%%%%%%%%%%%%%%%%%%%%%%%%%%%%%%%%%%%%%%%%%%%%%%%%%%%%%%%%%
%%%%%%%%%%%%%%%%%%%%%%%%%%%%%%%%%%%%%%%%%%%%%%%%%%%%%%%%%%%%%%%%%
\begin{corollary}
 Let  $A_1,\dots, A_n$ be AF-rings such that, for each $i =1,\dots, n$ and for each $M_i\in\Max(A_i)$, $\td(A_{iM_i}) = \td(A_{i})$. Then
 $$\dim(A_{1}\otimes\dots \otimes A_{n}) = \sum_{1\leq i\leq n}\td(A_{i})\ -\ \max\big\{\td(A_{i}) - \dim(A_{i})\big\}_{1\leq i\leq n}.$$
\end{corollary}

The above corollary recovers Wadsworth's aforementioned result.

%%%%%%%%%%%%%%%%%%%%%%%%%%%%%%%%%%%%%%%%%%%%%%%%%%%%%%%%%%%%%%%%%
%%%%%%%%%%%%%%%%%%%%%%%%%%%%%%%%%%%%%%%%%%%%%%%%%%%%%%%%%%%%%%%%%
\begin{corollary}[{\cite[Theorem 3.8]{W}}]
 Let  $D_1,\dots, D_n$ be AF-domains. Then
 $$\dim(D_{1}\otimes\dots \otimes D_{n}) = \sum_{1\leq i\leq n}\td(D_{i})\ -\ \max\big\{\td(D_{i}) - \dim(D_{i})\big\}_{1\leq i\leq n}.$$
\end{corollary}

Next, a sufficient condition involves the minimal prime ideals.

%%%%%%%%%%%%%%%%%%%%%%%%%%%%%%%%%%%%%%%%%%%%%%%%%%%%%%%%%%%%%%%%%
%%%%%%%%%%%%%%%%%%%%%%%%%%%%%%%%%%%%%%%%%%%%%%%%%%%%%%%%%%%%%%%%%
\begin{corollary}
 Let  $A_1,\dots, A_n$ be AF-rings such that, for each $i =1,\dots, n$ and for each minimal prime ideal $p_i$ of $A_i$, $\td(A_{i}/p_{i}) = \td(A_{i})$. Then
 $$\dim(A_{1}\otimes\dots \otimes A_{n}) = \sum_{1\leq i\leq n}\td(A_{i})\ -\ \max\big\{\td(A_{i}) - \dim(A_{i})\big\}_{1\leq i\leq n}.$$
\end{corollary}

%%%%%%%%%%%%%%%%%%%%%%%%%%%%%%%%%%%%%%%%%%%%%%%%%%%%%%%%%%%%%%%%%
%%%%%%%%%%%%%%%%%%%%%%%%%%%%%%%%%%%%%%%%%%%%%%%%%%%%%%%%%%%%%%%%%
\begin{corollary}
 Let  $A_1,\dots, A_n$ be equicodimensional AF-rings. Then
 $$\dim(A_{1}\otimes\dots \otimes A_{n}) = \sum_{1\leq i\leq n}\td(A_{i})\ -\ \max\big\{\td(A_{i}) - \dim(A_{i})\big\}_{1\leq i\leq n}.$$
\end{corollary}

It is known \cite[Corollary 3.3]{Giro} that if $A$ is an AF-ring, then
$$\dim(A\otimes A) = \dim_v(A\otimes A) \leq \dim(A) + \td(A) = \dim_v(A) +\td(A).$$
This result follows also from Corollary~\ref{af-cor1}. Now, applying Theorem~\ref{af-main2} to $A\otimes A$ we obtain:

%%%%%%%%%%%%%%%%%%%%%%%%%%%%%%%%%%%%%%%%%%%%%%%%%%%%%%%%%%%%%%%%%
%%%%%%%%%%%%%%%%%%%%%%%%%%%%%%%%%%%%%%%%%%%%%%%%%%%%%%%%%%%%%%%%%
\begin{corollary}\label{af-cor2}
Let  $A$  be an AF-ring. Then, the following assertions are equivalent:
\begin{enumerate}[label=(\rm\roman*)]
\item $\dim(A\otimes A) =\dim(A) + \td(A)$;
\item $\exists\ M,N\in\Max(A)$ with $\htt(M)=\dim(A)$, $\td(A_{N})=\td(A)$,  and $\td(A/N) \leq \td(A/M)$.
\end{enumerate}
\end{corollary}

Next, we provide an example of an AF-ring $A$ with $\dim(A\otimes A) \lneqq \dim(A) + \td(A)$.

%%%%%%%%%%%%%%%%%%%%%%%%%%%%%%%%%%%%%%%%%%%%%%%%%%%%%%%%%%%%%%%%%
%%%%%%%%%%%%%%%%%%%%%%%%%%%%%%%%%%%%%%%%%%%%%%%%%%%%%%%%%%%%%%%%%
\begin{example}
Let $K$ be a field extension of $k$ with $\td(K) = 2$ and let $A:=K\times k[X]$, where $X$ is an  indeterminate over $k$. Then, $A$ is an AF-ring with $\dim(A)=1$ and $\td(A)=2$. The maximal ideals of  $A$  are  $(0)\times k[X]$ and $K \times N$ with $N\in\Max(k[X])$.  Moreover, $\htt((0)\times k[X])= 0$; $\td(A_{(0)\times k[X]})= 2=\td(A)$;  $\htt(K\times N) = 1=\dim(A)$; $\td(A/((0)\times k[X]))=2$; and $\td(A/(K\times N))$=0.  So, by Corollary~\ref{af-cor2}, we obtain $\dim(A\otimes A)\lneqq \dim A + \td(A)= 3$.
\end{example}

We will conclude this subsection by an illustrative example which requires the following technical lemma.

%%%%%%%%%%%%%%%%%%%%%%%%%%%%%%%%%%%%%%%%%%%%%%%%%%%%%%%%%%%%%%%%%
%%%%%%%%%%%%%%%%%%%%%%%%%%%%%%%%%%%%%%%%%%%%%%%%%%%%%%%%%%%%%%%%%
\begin{lemma}\label{af-lem2}
Let $A$ be an AF-ring such that there exist $p,q\in\Spec(A)$ with $\td(A_{p})\not= \td(A_{q})$. Then, for any AF-ring  $B$, $A\otimes B$  is not the tensor product of a finite number of AF-domains.
\end{lemma}

%%%%%%%%%%%%%%%%%%%%%%%%%%%%%%%%%%%%%%%%%%%%%%%%%%%%%%%%%%%%%%%%%
%%%%%%%%%%%%%%%%%%%%%%%%%%%%%%%%%%%%%%%%%%%%%%%%%%%%%%%%%%%%%%%%%
\begin{example}
For each integer $n\geq1$, there exist two AF-rings $A_1$ and $A_2$ such that:
\begin{enumerate}
\item $\dim(A_1\otimes A_2) = n$;
\item $A_1\otimes A_2$  is not the tensor product of a finite number of AF-domains;
\item If, in addition, there exists a non-finitely generated separable extension of $k$, then neither $A_1$ nor $A_2$ is a finite
direct product of AF-domains.
\end{enumerate}

Indeed, let $K$ be a separable extension of $k$. Let $V_1 := K(X)[Y]_{(Y)} = K(X) + M_1$, $V:=K(Y)[X]_{(X)}= K(Y) + M$, and $V_2 := K[Y]_{(Y)}+ M= K + M_2$. Then,  $V_1$ and $V_2$ are, respectively, one-dimensional and two-dimensional valuation domains of $K(X,Y)$. Since  $V_1$  and $V_2$  are incomparable,  $T := V_1\cap V_2$ is a two-dimensional Pr\"ufer domain with only two maximal ideals, $M_1$  and $M_2$, such that $T_{M_1}= V_1$ and $T_{M_2}= V_2$ \cite[Theorem 11.11]{Na}. Let $I := M_1M_2$ and $R := T/I$. Then, $R$ is a  zero-dimensional ring (and, a fortiori, an AF-ring) with  only two prime ideals, $p_1:= M_1/I$ and $p_2:= M_2/I$. Further, $\td(R/p_1) = 1$  and $\td(R/p_2) = 0$. By Corollary~\ref{af-cor2}, we have
$$\dim(R\otimes R[n]) =\dim((R\otimes R)[n]) =\dim(R\otimes R) + n = \td(R) + n = 1 +n.$$
Moreover, by Lemma~\ref{af-lem2}, $R\otimes R[n]$  is not the tensor product of a finite number of $AF$-domains; so it suffices to take $A_1 := R$ and $A_2 := R[n-1]$.

Now assume that $K$ is not finitely generated over $k$. So  $K\otimes K$ is a reduced \cite[Theorem 39]{ZS}, zero-dimensional \cite[Theorem 3.1]{Sh}, and
non-Noetherian \cite[Theorem 11]{Va} ring. Then $\Spec(K\otimes K)$ is infinite \cite[Lemma 0]{Va}. Next, let $A := K\otimes R$. Since $A$ is an integral extension of $R$, it is zero-dimensional. Moreover, there exist two prime ideals of $A$, $P_1$ and $P_2$ such that $P_1\cap R = p_1$ and $P_2\cap R = p_2$ with $\td(A/P_1)=1$ and $\td(A/P_2) = 0$. Since $K$ is the quotient field of  $R/p_2$ and $\Spec(K\otimes K)$ is infinite, by \cite[Proposition 3.2]{W}, $\Spec(A)$ is infinite. So $A$ is not a finite direct product of AF-domains and the  same holds for $A[n]$. By Corollary~\ref{af-cor2}, $\dim(A\otimes A[n]) =\dim((A\otimes A)[n]) =\dim(A\otimes A) + n = \td(A) + n = 1 +n$. Moreover, by Lemma~\ref{af-lem2},
$A\otimes A[n]$ is not the tensor product  of a finite number of AF-domains. So it  suffices to take $A_1 := A$  and $A_2 := A[n-1]$, completing the proof of the example.
\end{example}

%%%%%%%%%%%%%%%%%%%%%%%%%%%%%%%%%%%%%%%%%%%%%%%%%%%%%%%%%%%%%%%%%%%%%%%%%%%%%%%%%%%%%%%%%%%%%%%%%%%%%
%%%%%%%%%%%%%%%%%%%%%%%%%%%%%%%%%%%%%%%%%%%%%%%%%%%%%%%%%%%%%%%%%%%%%%%%%%%%%%%%%%%%%%%%%%%%%%%%%%%%%
%%%%%%%%%%%%%%%%%%%%%%%%%%%%%%%%%%%%%%%%%%%%%%%%%%%%%%%%%%%%%%%%%%%%%%%%%%%%%%%%%%%%%%%%%%%%%%%%%%%%%
%%%%%%%%%%%%%%%%%%%%%%%%%%%%%%%%%%%%%%%%%%%%%%%%%%%%%%%%%%%%%%%%%%%%%%%%%%%%%%%%%%%%%%%%%%%%%%%%%%%%%
\subsection{Transfer of the (locally) Jaffard property}

\noindent In this subsection, we first establish the transfer of the locally Jaffard property in some context of tensor products. Then, we give some  formulas for computing the valuative dimension of the tensor product of an AF-ring and an arbitrary ring. We conclude with the fact that the tensor product of an AF-ring and a Jaffard ring is not necessarily a Jaffard ring.

Next, we announce the main result of this subsection. Notice that the tensor product of two AF-rings is locally Jaffard (since it is an AF-ring).

%%%%%%%%%%%%%%%%%%%%%%%%%%%%%%%%%%%%%%%%%%%%%%%%%%%%%%%%%%%%%%%%%
%%%%%%%%%%%%%%%%%%%%%%%%%%%%%%%%%%%%%%%%%%%%%%%%%%%%%%%%%%%%%%%%%
\begin{theorem}\label{af-main3}
Let $A$ be an AF-ring and $B$ a locally Jaffard ring. Then, $A\otimes B$ is a locally Jaffard ring.
\end{theorem}

The proof of this result lies on a very important lemma which correlates the height of a prime ideal  of $A\otimes B$ to its traces on $A$ and $B$ via the transcendence degrees; namely, we have:

%%%%%%%%%%%%%%%%%%%%%%%%%%%%%%%%%%%%%%%%%%%%%%%%%%%%%%%%%%%%%%%%%
%%%%%%%%%%%%%%%%%%%%%%%%%%%%%%%%%%%%%%%%%%%%%%%%%%%%%%%%%%%%%%%%%
\begin{lemma}
Let $A$ be an AF-ring and $B$ an arbitrary ring. Let $P\in\Spec(A\otimes B)$ with $p := P\cap A$ and $q := P\cap B$. Then
$$\htt(P) + \td\big((A\otimes B)/P\big) =\td(A_p) + \htt(q[\td(A_p)]) + \td(B/q).$$
\end{lemma}

Next, we give some applications of Theorem~\ref{af-main3}. The first one establishes a formula for the valuative dimension of $A\otimes B$ where $A$ is an AF-ring. To this purpose, one should first examine the variation of the function $D$ between $B$ and its associated polynomial rings.

%%%%%%%%%%%%%%%%%%%%%%%%%%%%%%%%%%%%%%%%%%%%%%%%%%%%%%%%%%%%%%%%%
%%%%%%%%%%%%%%%%%%%%%%%%%%%%%%%%%%%%%%%%%%%%%%%%%%%%%%%%%%%%%%%%%
\begin{lemma}
Let $A$ be an AF-ring, $B$ an arbitrary ring, and $p\in\Spec(A)$. Then, for any $n\geq 1$, we have
$$D\big(\td(A_{p}), \htt(p), B[n]\big) = D\big(\td(A_{p}) + n,  \htt(p) + n,  B\big).$$
\end{lemma}

The next result provides a formula for the valuative dimension.

%%%%%%%%%%%%%%%%%%%%%%%%%%%%%%%%%%%%%%%%%%%%%%%%%%%%%%%%%%%%%%%%%
%%%%%%%%%%%%%%%%%%%%%%%%%%%%%%%%%%%%%%%%%%%%%%%%%%%%%%%%%%%%%%%%%
\begin{corollary}\label{af-cor3}
Let $A$ be an AF-ring and $B$ a ring with finite valuative dimension $\geq 1$. Then, for any $r\geq \dim_{v}(B)-1$, we have
$$\dim_{v}(A\otimes B)= \max\big\{\htt(q[r])+\min\big(\td(A_{p}),\htt(p)+\td(B/q)\big)\big\}$$
where $(p,q)$ ranges over $\Spec(A)\times\Spec(B)$.
\end{corollary}

The special case where $A$ is an AF-domain  yields a more simplified formula.

%%%%%%%%%%%%%%%%%%%%%%%%%%%%%%%%%%%%%%%%%%%%%%%%%%%%%%%%%%%%%%%%%
%%%%%%%%%%%%%%%%%%%%%%%%%%%%%%%%%%%%%%%%%%%%%%%%%%%%%%%%%%%%%%%%%
\begin{corollary}
Let $A$ be an AF-domain and $B$ a ring with finite valuative dimension $\geq 1$. Then, for any $r\geq \dim_{v}(B)-1$, we have
$$\dim_{v}(A\otimes B)= \max\big\{\htt(q[r])+\min\big(\td(A),\dim(A)+\td(B/q)\big)\big\}$$
where $q$ ranges over $\Spec(B)$.
\end{corollary}

The next two results feature special contexts where the tensor product is Jaffard.

%%%%%%%%%%%%%%%%%%%%%%%%%%%%%%%%%%%%%%%%%%%%%%%%%%%%%%%%%%%%%%%%%
%%%%%%%%%%%%%%%%%%%%%%%%%%%%%%%%%%%%%%%%%%%%%%%%%%%%%%%%%%%%%%%%%
\begin{corollary}
Let $A$ be an AF-domain and $B$ a ring such that $\dim_{v}(B)\leq \td(A)+1$. Then, $A\otimes B$ is a Jaffard ring.
\end{corollary}

 Recall that, for any ring $B$ of valuative dimension 2, the ring $B[X]$ is locally Jaffard \cite[Proposition 1(ii)]{Cahen}. Also, \cite[Example 3.2]{ABDFK} is an
example of a Jaffard ring $B$ that is not locally  Jaffard but $B[X]$ is locally Jaffard.

%%%%%%%%%%%%%%%%%%%%%%%%%%%%%%%%%%%%%%%%%%%%%%%%%%%%%%%%%%%%%%%%%
%%%%%%%%%%%%%%%%%%%%%%%%%%%%%%%%%%%%%%%%%%%%%%%%%%%%%%%%%%%%%%%%%
\begin{corollary}
Let $A$ be an AF-domain and $B$ a Jaffard ring such that $B[X]$ is locally Jaffard. Then, $A\otimes B$ is a Jaffard ring.
\end{corollary}

We close this section with an example where the tensor product of  an AF-domain and a Jaffard ring is not necessarily a Jaffard ring. This shows that a similar result to Theorem~\ref{af-main3} does not hold, in general, for the transfer of the Jaffard property.
%%%%%%%%%%%%%%%%%%%%%%%%%%%%%%%%%%%%%%%%%%%%%%%%%%%%%%%%%%%%%%%%%
%%%%%%%%%%%%%%%%%%%%%%%%%%%%%%%%%%%%%%%%%%%%%%%%%%%%%%%%%%%%%%%%%
\begin{example}
We deduce this example from \cite{ABDFK}. Let $Z_1,Z_2,Z_3,Z_4$ be four  indeterminates over  $k$. Let $L:= k(Z_1,Z_2,Z_3,Z_4)$ and $V_1 := k(Z_1,Z_2,Z_3)[Z_4]_{(Z_4)}$. Then, $V_1$ is a one-dimensional valuation ring of $L$ with maximal ideal $M_1 := Z_4V_1$. Let $V'$ be a one-dimensional valuation overring of $k(Z_4)[Z_2,Z_3]$ of the form
$V' := k(Z_4) + M'$ and  $V{'}_2 = k[Z_4]_{(Z_4)}+ M' = k + M'_2$, where $M'_2 = Z_4k[Z_4]_{(Z_4)} + M'$. So, $V'_2$ is a two- dimensional valuation ring. Now, let  $V = k(Z_2,Z_3,Z_4)[Z_1]_{(Z_1)} = k(Z_2,Z_3,Z_4) + M$, with $M= Z_1V$, and $M_2 = M'_2 + M$. Then, $V_2 := V'_2 + M =k + M_2$ is a three-dimensional valuation ring.

We claim that $V_1$ and $V_2$ are incomparable. Otherwise, $V_2\subset V_1$ and hence $V_1 = (V_2)_M$. So, $M$ is a divided prime ideal of $V_2$. That is, $Z_4V_1 = M_1 = MV_{2M}$. So $1 = Z_4 Z_4^{-1}\in MV = M$, the desired contradiction. Now, $V_1$  and $V_2$  have the same quotient field $L$. By \cite[Theorem 11.11]{Na}, $S := V_1\cap V_2$ is a three-dimensional Pr\"ufer domain with only two maximal ideals, $N_1$ and $N_2$, such that $S_{N_1}= V_1$ and $S_{N_2}= V_2$. Next, let $F := k(Z_1)$, $f : V_1\longrightarrow
k(Z_1,Z_2,Z_3)$ be the natural ring homomorphism, and $D :=f^{-1}(F) = F + M_1$. Let $g :S\longrightarrow S/N_1\cong V_1/N_1\cong k(Z_1,Z_2,Z_3)$ be the natural ring homomorphism  and $B:=g^{-1}(F)$. We have $B = D\cap S = D\cap V_2$  and $\dim(B)
=\dim(S) = 3$. Moreover, by  \cite[Theorem 2.11]{ABDFK}, we obtain
$$\dim_v(B) =\max\big\{\dim_v(S),\dim_v(F) + \dim_v(S_{N_1}) + \t.d.(S /N_1:F)\big\} = 3.$$
Therefore, $B$ is Jaffard. Since $B = D\cap V_2$ and $V_1$, $V_2$  are incomparable, it follows that $B_{n_1} = D$  and $B_{n_2}= V_2$, where $\{n_1, n_2\}=\Max(B)$. Moreover,
$\htt(n_1[s]) = \htt(n_1B_{n_1}[s]) = \htt(M_1[s])$, for any positive integer $s$. Since $V_1$ is Jaffard, by \cite[Theorem 1.7]{Ay}, $\htt  _{D[s]}(M_1[s])=\htt _{V_1}(M_1)+\min(s,2)$.
Then $\htt(n_1) = 1$, $\htt(n_1[X_1]) = 2$, and  $\htt(n_1[X_1,X_2]) = 3$;
$\td(B/n_1) = \td(D/M_1) = 1$, and $\td(B/n_2)  = \td(V_2/M_2) = 0$.

Let $A := k(X)$. By Theorem~\ref{af-main1}, we have
$$\dim(A\otimes B) = D(\td(A), 0, B) =\max\big\{\htt(q[X_1]) +\min\big(1, \td(B/q)\big) \mid q\in\Spec(B)\big\}.$$
For $q:=  n_1$, it yields  $\htt(n_1[X_1]) + \min(1,\td(B/n_1))= 2 + 1 = 3$; for $q := n_2$, it yields $\htt(n_2[X_1]) +\min(1,\td(B/n_2)) = \htt(n_2) = 3$,  and  $\htt(q[X_1]) +\min(1,\td(B/q))\leq 3$  for every prime ideal $q$ of $B$ contained in $n_2$. Consequently, $\dim(A\otimes B) = 3$. By Corollary~\ref{af-cor3}, $\dim_v(A\otimes B) =\max\big\{\htt(q[X_1,X_2]) +\min(1, \td(B/q))\mid q\in\Spec(B)\big\}$. For $q := n_1$, it is $\htt(n_1[X_1,X_2])+ \min(1, \td(B/n_1)) = 3 + 1 = 4$. Therefore, $\dim_v(A\otimes B) = 4\not= \dim(A\otimes B)$. Consequently, $A\otimes B$ is not a Jaffard ring, completing the proof of the example.
\end{example}

%\newpage
%%%%%%%%%%%%%%%%%%%%%%%%%%%%%%%%%%%%%%%%%%%%%%%%%%%%%%%%%%%%%%%%%%%%%%%%%%%%%%%%%%%%%%%%%%%%%%%%%%%%%%%%%%%%%%%%%%%%%%%%%%%%%%%%%%%%%%%%%%%%%%
%%%%%%%%%%%%%%%%%%%%%%%%%%%%%%%%%%%%%%%%%%%%%%%%%%%%%%%%%%%%%%%%%%%%%%%%%%%%%%%%%%%%%%%%%%%%%%%%%%%%%%%%%%%%%%%%%%%%%%%%%%%%%%%%%%%%%%%%%%%%%%
%%%%%%%%%%%%%%%%%%%%%%%%%%%%%%%%%%%%%%%%%%%%%%%%%%%%%%%%%%%%%%%%%%%%%%%%%%%%%%%%%%%%%%%%%%%%%%%%%%%%%%%%%%%%%%%%%%%%%%%%%%%%%%%%%%%%%%%%%%%%%%
%%%%%%%%%%%%%%%%%%%%%%%%%%%%%%%%%%%%%%%%%%%%%%%%%%%%%%%%%%%%%%%%%%%%%%%%%%%%%%%%%%%%%%%%%%%%%%%%%%%%%%%%%%%%%%%%%%%%%%%%%%%%%%%%%%%%%%%%%%%%%%
\section{Tensor products of AF-rings over a zero-dimensional ring}\label{z}

\noindent  This section is devoted to \cite{BGK2}. Its purpose is to extend all the known results on the dimension of tensor products of AF-rings over a field  to the general case of AF-rings over a zero-dimensional ring.

Throughout this section, $R$ denotes a zero-dimensional ring, and algebras (resp., tensor products), when not specifically marked, are taken over (resp., relative to) $R$. We denote by $(A,\lambda _A)$ an algebra $A$ and its associated ring homomorphism $\lambda _A:R\rightarrow A$; and, by $\lambda _A^*$, the associated spectral map $\Spec(A)\rightarrow \Spec(R)$. Notice that for any prime ideal $P$ of $A$, $\lambda _A^{-1}(P)$ is a maximal ideal of $R$. So, we define the transcendence degree of the algebra $A$ over $R$ as follows
$$\t.d.(A:_{\lambda _A}R)=\max\big\{\t.d.\big(A/P:R/\lambda _A^{-1}(P)\big)\mid P\in\Spec(A)\big\}.$$
We write $\td(A:R)$ or just $\td(A)$  as an abbreviation for $\t.d.(A:_{\lambda_A}R)$, when there is no ambiguity. All along this section, we consider only
algebras $(A,\lambda _A)$ such that $\td(A)<\infty$, which also ensures that $\dim(A)<\infty$. If $A$ is an integral domain, $p_A$ denotes $\Ker(\lambda_A)$.

First of all, observe that the transcendence degree of an algebra $A$ depends on its $R$-module structure, as shown by the next example.

%%%%%%%%%%%%%%%%%%%%%%%%%%%%%%%%%%%%%%%%%%%%%%%%%%%%%%%%%%%%%%%%%
%%%%%%%%%%%%%%%%%%%%%%%%%%%%%%%%%%%%%%%%%%%%%%%%%%%%%%%%%%%%%%%%%
\begin{example}
Let $R:=k(X)\times k$ and $A:=k(X)$, where $k$ is a field. Consider the two ring homomorphisms $\lambda _1:R\rightarrow A$ and $\lambda _2: R\rightarrow A$ defined by $\lambda_1(x,y)=x$ and $\lambda_2(x,y)=y$. Then $\t.d.(A:_{\lambda_1}R)=\t.d.(k(X):k(X))=0$ whereas $\t.d.(A:_{\lambda_2}R)=\t.d.(k(X):k)=1$.
\end{example}

The following lemma provides simple generalizations of well-known facts for algebras over a field.

%%%%%%%%%%%%%%%%%%%%%%%%%%%%%%%%%%%%%%%%%%%%%%%%%%%%%%%%%%%%%%%%%
%%%%%%%%%%%%%%%%%%%%%%%%%%%%%%%%%%%%%%%%%%%%%%%%%%%%%%%%%%%%%%%%%
\begin{lemma}
Let $(A,\lambda _A)$ be an algebra, $P\in\Spec(A)$, and $p:=\lambda _A^{-1}(P)$. Then:
\begin{enumerate}
\item $\htt(P)+\td(A/P:R)\leq \td(A_P:R)=\td((A/pA)_{P/pA}:R)$.
\item $\htt(P)=\htt(P/pA)$.
\item $\htt(P[n])=\htt((P/pA)[n])$, for each $n\geq 1$.
\item If $A$ is locally Jaffard, then so is $A/qA$, for each $q\in\Spec(R)$ with $qA \not=A$.
\end{enumerate}
\end{lemma}

%%%%%%%%%%%%%%%%%%%%%%%%%%%%%%%%%%%%%%%%%%%%%%%%%%%%%%%%%%%%%%%%%%%%%%%%%%%%%%%%%%%%%%%%%%%%%%%%%%%%%
%%%%%%%%%%%%%%%%%%%%%%%%%%%%%%%%%%%%%%%%%%%%%%%%%%%%%%%%%%%%%%%%%%%%%%%%%%%%%%%%%%%%%%%%%%%%%%%%%%%%%
%%%%%%%%%%%%%%%%%%%%%%%%%%%%%%%%%%%%%%%%%%%%%%%%%%%%%%%%%%%%%%%%%%%%%%%%%%%%%%%%%%%%%%%%%%%%%%%%%%%%%
%%%%%%%%%%%%%%%%%%%%%%%%%%%%%%%%%%%%%%%%%%%%%%%%%%%%%%%%%%%%%%%%%%%%%%%%%%%%%%%%%%%%%%%%%%%%%%%%%%%%%
\subsection{Tensorially compatible algebras}

\noindent Let $(A_1,\lambda _1)$ and $(A_2,\lambda _2)$ be algebras. For $i=1,2,$ we denote by $\mu _i:A_i\rightarrow A_1\otimes   A_2$ the canonical $A_i$-algebra homomorphism. The algebra  $A_1\otimes   A_2$, when not specifically indicated, has $\lambda_{A_1\otimes A_2} = \mu _1\circ \lambda_1=\mu _2\circ \lambda_2$ as its associated ring homomorphism.  Finally, let
$$\Gamma(A_1,A_2):=\big\{(P_1,P_2) \in\Spec (A_1)\times\Spec (A_2)\mid  \lambda_1^{-1}(P_1)=\lambda_2^{-1}(P_2)\big\}.$$

We are interested in algebras $(A_1,\lambda _1)$ and $(A_2,\lambda _2)$ such that $A_1\otimes A_2\neq 0$, and call such algebras \emph{tensorially compatible}. The next result provides some elementary and useful characterizations of tensorially compatible algebras. For a more general result, we refer the reader to \cite[Corollary 3.2.7.1]{GD}.

%%%%%%%%%%%%%%%%%%%%%%%%%%%%%%%%%%%%%%%%%%%%%%%%%%%%%%%%%%%%%%%%%
%%%%%%%%%%%%%%%%%%%%%%%%%%%%%%%%%%%%%%%%%%%%%%%%%%%%%%%%%%%%%%%%%
\begin{proposition}
Let $(A_1,\lambda _1)$ and $(A_2,\lambda _2)$ be algebras. Then, the following
assertions are equivalent:
\begin{enumerate}
\item $(A_1,\lambda _1)$ and $(A_2,\lambda _2)$ are tensorially compatible;
\item $\lambda _1^*(\Spec (A_1))\cap \lambda _2^*(\Spec (A_2))\neq\emptyset$;
\item $\exists\ P_1\Spec(A_1)$ such that $\lambda_1^{-1}(P_1)A_2\neq A_2$;
\item $\exists\ P_2\Spec(A_2)$ such that $\lambda_2^{-1}(P_2)A_1\neq A_1$;
\item $\exists\ p\Spec(R)$ such that $pA_1\neq A_1$ and $pA_2\neq A_2$;
\item $\Ker(\lambda_1)+\Ker(\lambda_2)\neq R$.
\end{enumerate}
\end{proposition}

A similar result holds for any finite number of algebras, as shown below.

%%%%%%%%%%%%%%%%%%%%%%%%%%%%%%%%%%%%%%%%%%%%%%%%%%%%%%%%%%%%%%%%%
%%%%%%%%%%%%%%%%%%%%%%%%%%%%%%%%%%%%%%%%%%%%%%%%%%%%%%%%%%%%%%%%%
\begin{proposition}
Let $(A_1,\lambda_1),\dots ,(A_n,\lambda_n)$ be algebras. Then, the following assertions are equivalent:
\begin{enumerate}
\item $A_1\otimes \cdots \otimes A_n\neq \boldsymbol {0}$;
\item $\lambda _1^*(\Spec (A_1))\cap \lambda _2^*(\Spec (A_2))\cap \cdots
\cap \lambda _n^*(\Spec (A_n))\neq \emptyset$;
\item $\exists\ p\Spec(R)$ such that $pA_i\neq A_i$, for each $i=1,2,\dots,n$.
\end{enumerate}
\end{proposition}

The next result establishes an analogue to \cite[Proposition 2.3]{W}.

%%%%%%%%%%%%%%%%%%%%%%%%%%%%%%%%%%%%%%%%%%%%%%%%%%%%%%%%%%%%%%%%%
%%%%%%%%%%%%%%%%%%%%%%%%%%%%%%%%%%%%%%%%%%%%%%%%%%%%%%%%%%%%%%%%%
\begin{proposition}
Let $(A_1,\lambda_1)$ and $(A_2,\lambda_2)$ be algebras and $(P_1,P_2)\in \Spec(A_1)\times\Spec(A_2)$  with $\lambda_1^{-1}(P_1)=\lambda_2^{-1}(P_2)=p$. Let
$\Omega :=\big\{Q\in \Spec(A_1\otimes A_2)\mid\mu_i^{-1}(Q)=P_i,\ i=1,2\big\}$.
Then
\begin{enumerate}
\item $\Omega$ is lattice isomorphic to $\Spec\big(\frac{A_{P_1}}{P{_1}A_{P_1}}\otimes_{\frac{R}{p}}\frac{A_{P_2}}{P{_2}A_{P_2}}\big)$.
\item $Q\in\Omega$ is minimal in $\Omega$ if and only if $\td((A_1\otimes A_2)/Q)=\td(A_1/P_1)+\td(A_2/P_2)$.
\item If $Q_o\in \Spec (A_1\otimes A_2)$ and $\mu_i^{-1}(Q_o)\supseteq P_i$ ($i=1,2$),
then $\exists\ Q\in\Omega$ such that $Q\subseteq Q_o$.
\end{enumerate}
\end{proposition}

Follow two applications of the above result, which extend two known results on algebras over a field \cite{W} to $R$-algebras.

%%%%%%%%%%%%%%%%%%%%%%%%%%%%%%%%%%%%%%%%%%%%%%%%%%%%%%%%%%%%%%%%%
%%%%%%%%%%%%%%%%%%%%%%%%%%%%%%%%%%%%%%%%%%%%%%%%%%%%%%%%%%%%%%%%%
\begin{corollary}
Let $(A_1,\lambda _1)$ and $(A_2,\lambda _2)$ be tensorially compatible algebras and
let $Q\in\Spec(A_1\otimes A_2)$. Then
$$\htt(Q)\geq \htt(\mu_1^{-1}(Q))+\htt(\mu_2^{-1}(Q)).$$
\end{corollary}

%%%%%%%%%%%%%%%%%%%%%%%%%%%%%%%%%%%%%%%%%%%%%%%%%%%%%%%%%%%%%%%%%
%%%%%%%%%%%%%%%%%%%%%%%%%%%%%%%%%%%%%%%%%%%%%%%%%%%%%%%%%%%%%%%%%
\begin{corollary}\label{z-cor1}
Let $(A_1,\lambda _1)$ and $(A_2,\lambda _2)$ be tensorially compatible algebras.
Then
$$\td(A_1\otimes A_2) =\max\big\{\td(A_1/P_1)+\td(A_2/P_2)\mid (P_1,P_2)\in\Gamma\big\}\leq \td(A_1)+\td(A_2).$$
\end{corollary}

Let $(A_1,\lambda _1)$ and $(A_2,\lambda _2)$ be tensorially compatible algebras. Clearly,
$\td(A_1\otimes A_2)=\td(A_1)+\td(A_2)$ if and only if there exists $(P_1,P_2)\in\Gamma$ with $\td(A_1)=\td(A_1/P_1)$ and $\td(A_2)=\td(A_2/P_2)$. The second condition holds, for instance, if $A_1$ and $A_2$ are integral domains or if $\Spec(R)$ is reduced to only one prime ideal. In general, the equality fails as it is shown by the next example. Moreover, when $R$ is a field, we have $\dim(A_1\otimes A_2)\geq \dim(A_1)+\dim(A_2)$ \cite[Corollary 2.5]{W}. This is not true, in general, in the zero-dimensional case, as shown below.

%%%%%%%%%%%%%%%%%%%%%%%%%%%%%%%%%%%%%%%%%%%%%%%%%%%%%%%%%%%%%%%%%
%%%%%%%%%%%%%%%%%%%%%%%%%%%%%%%%%%%%%%%%%%%%%%%%%%%%%%%%%%%%%%%%%
\begin{example}
Let $R:=\R\times\R$, $A_1:=\R$ and $A_2:=\R\times\R[X]$. Consider the two ring homomorphisms $\lambda_1:R\rightarrow A_1$ and $\lambda_2:R\rightarrow A_2$ defined by $\lambda_1(x,y)=x$ and $\lambda_2(x,y)=(x,y)$. Clearly, $A_1$ and $A_2$ are tensorially compatible. We claim that
$$\td(A_1\otimes A_2) \lneqq \td(A_1)+\td(A_2)\
\text{ and }\
\dim(A_1\otimes A_2)\lneqq \dim(A_1)+\dim(A_2).$$

Indeed, one can easily see that $\td(A_1)=\td(\R:_{\lambda _1}R)=\t.d.(\R:\R)=0$, and $\td (A_2)= \td(\R\times \R[X]:_{\lambda _2}R)=\max\big\{\t.d.(\R:\R), \t.d.(\R[X]:\R)\big\}=1$. Moreover, by Corollary~\ref{z-cor1}, $\td(A_1\otimes A_2) = \td(A_1) + \td\big(A_2/(0\times\R[X])\big)=\t.d.(\R:\R)+\t.d.(\R:\R)=0$. It follows that
$$\dim(A_1\otimes A_2)\leq \td(A_1\otimes A_2)=0\lneqq \td(A_1)+\td(A_2)=\dim(A_1)+\dim(A_2)=1$$
completing the proof of the example.
\end{example}

%%%%%%%%%%%%%%%%%%%%%%%%%%%%%%%%%%%%%%%%%%%%%%%%%%%%%%%%%%%%%%%%%%%%%%%%%%%%%%%%%%%%%%%%%%%%%%%%%%%%%
%%%%%%%%%%%%%%%%%%%%%%%%%%%%%%%%%%%%%%%%%%%%%%%%%%%%%%%%%%%%%%%%%%%%%%%%%%%%%%%%%%%%%%%%%%%%%%%%%%%%%
%%%%%%%%%%%%%%%%%%%%%%%%%%%%%%%%%%%%%%%%%%%%%%%%%%%%%%%%%%%%%%%%%%%%%%%%%%%%%%%%%%%%%%%%%%%%%%%%%%%%%
%%%%%%%%%%%%%%%%%%%%%%%%%%%%%%%%%%%%%%%%%%%%%%%%%%%%%%%%%%%%%%%%%%%%%%%%%%%%%%%%%%%%%%%%%%%%%%%%%%%%%
\subsection{Krull dimension}

\noindent This subsection investigates the Krull dimension of tensor products of AF-rings over zero-dimensional rings. We first extend Wadsworth's definition of AF-rings over fields to AF-rings over zero-dimensional rings. Recall that $R$ denotes a zero-dimensional ring and algebras are taken over $R$.

%%%%%%%%%%%%%%%%%%%%%%%%%%%%%%%%%%%%%%%%%%%%%%%%%%%%%%%%%%%%%%%%%
%%%%%%%%%%%%%%%%%%%%%%%%%%%%%%%%%%%%%%%%%%%%%%%%%%%%%%%%%%%%%%%%%
\begin{definition}
Under the above notation, an algebra $(A,\lambda_A)$ is an AF-ring if
$$\htt(P)+\td(A/P)=\td(A_P),\ \forall\ P\in\Spec(A).$$
\end{definition}

It is worthwhile observing that this notion of AF-ring is independent of the structure of algebra defined by the ring homomorphism $\lambda_A$. Indeed, let $A$ be an algebra and let $\lambda $ and  $\lambda' $ be two ring homomorphisms defining two different structures of algebra over $R$ on $A$. Let $P\in\Spec(A)$ and $\pi:A\rightarrow A/P$ be the canonical
ring homomorphism. Let $p:=\Ker(\pi\circ\lambda)=\lambda^{-1}(P)$ and
$q:=\Ker(\pi\circ\lambda')=\lambda'^{-1}(P)$. One can regard $R/p$ and $R/q$ as subfields of $A/P$. Let $k:=R/p \cap R/q$. On one hand, we have $$\td(A/P:_{\lambda}R)=\t.d.(A/P:R/p)=\t.d.(A/P:k)-\t.d.(R/p:k)$$
and
$$\td(A/P:_{\lambda'}R)=\t.d.(A/P:R/q)=\t.d.(A/P:k)-\t.d.(R/q:k).$$
On the other hand, we have
\begin{eqnarray*}
\td(A_P:_\lambda R) &=  &\max\big\{\t.d.(A/Q:R/p) \mid Q\in \Spec(A)\ \text{and}\ Q\subseteq P\big\}\\
                    &=  &\max\big\{\t.d.(A/Q:k) \mid Q\in \Spec(A)\ \text{and}\ Q\subseteq P\big\} - \t.d.(R/p:k)
\end{eqnarray*}
and
\begin{eqnarray*}
\td(A_P:_{\lambda'} R) &=  &\max\big\{\t.d.(A/Q:R/q) \mid Q\in \Spec(A)\ \text{and}\ Q\subseteq P\big\}\\
                    &=  &\max\big\{\t.d.(A/Q:k) \mid Q\in \Spec(A)\ \text{and}\ Q\subseteq P\big\} - \t.d.(R/q:k).
\end{eqnarray*}
It follows that $\td(A_P:_\lambda R)-\td(A/P:_{\lambda}R)=\td(A_P:_{\lambda'} R)-\td(A/P:_{\lambda'}R)$. That is, $(A,\lambda)$ is an AF-ring if and only if $(A,\lambda')$ is an AF-ring.
\bigskip

Next, we provide some examples and basic properties of AF-rings.

%%%%%%%%%%%%%%%%%%%%%%%%%%%%%%%%%%%%%%%%%%%%%%%%%%%%%%%%%%%%%%%%%
%%%%%%%%%%%%%%%%%%%%%%%%%%%%%%%%%%%%%%%%%%%%%%%%%%%%%%%%%%%%%%%%%
\begin{lemma}
Let $\mathcal{R}$ be the class of AF-rings (over $R$) and let $n$ be a positive integer. Then:
\begin{enumerate}
\item $A\in\mathcal{R}$ $\Leftrightarrow$ $A/pA$ is an AF-ring over the field $R/p$, $\forall\ p\in\Spec(R)$ with $pA\neq A$.
\item Any finitely generated $R$-algebras and its integral extensions belong to $\mathcal{R}$.
\item If $A\in\mathcal{R}$, then $S^{-1}A\in\mathcal{R}$, for every multiplicative subset $S$ of $A$.
\item If $A\in\mathcal{R}$, then $A[n]\in\mathcal{R}$.
\item If  $A_1,\dots, A_n\in\mathcal{R}$ and are tensorially compatible, then $A_1\otimes\dots\otimes A_n\in\mathcal{R}$.
\item If  $A_1,\dots, A_n\in\mathcal{R}$, then $A_1\times\dots\times A_n\in\mathcal{R}$.
\item If $A\in\mathcal{R}$, then $A$ is locally Jaffard.
\end{enumerate}
\end{lemma}

Next, we establish adequate analogues of the main  results stated in Section~\ref{af} on the dimension of tensor products of AF-rings over a field. The first result provides a formula for the Krull dimension of the tensor product $A\otimes B$, where $A$ is an AF-ring.

%%%%%%%%%%%%%%%%%%%%%%%%%%%%%%%%%%%%%%%%%%%%%%%%%%%%%%%%%%%%%%%%%
%%%%%%%%%%%%%%%%%%%%%%%%%%%%%%%%%%%%%%%%%%%%%%%%%%%%%%%%%%%%%%%%%
\begin{theorem}\label{z-main1}
Let $A$ be an AF-ring and $B$ an algebra with $A\otimes B\neq 0$. Then
$$\dim(A\otimes B)=\max\big\{\htt(Q[\td(A_{P})])+\min\big(\td(A_{P}), \htt(P)+\td(B/Q)\big)\mid (P,Q)\in\Gamma(A,B)\big\}.$$
\end{theorem}

It is worthwhile noting that $\dim(A\otimes B)$ depends on the $R$-module structure of $A$ and $B$, as shown by the next example.

%%%%%%%%%%%%%%%%%%%%%%%%%%%%%%%%%%%%%%%%%%%%%%%%%%%%%%%%%%%%%%%%%
%%%%%%%%%%%%%%%%%%%%%%%%%%%%%%%%%%%%%%%%%%%%%%%%%%%%%%%%%%%%%%%%%
\begin{example}
Let $(A,\lambda_A)$ be an AF-ring and $(B,\lambda_B)$ an algebra with $\dim(A\otimes B)\neq 0$. Let $p\in\Spec(R)$ and $\pi:R\rightarrow R/p$ be the canonical ring homomorphism. On one hand, let $\lambda_1:R\times R\times R\rightarrow R/p\times A$ and $\lambda_2:R\times R\times R\rightarrow R/p\times B$ be ring homomorphisms defined by  $\lambda_1(x,y,z)= (\pi(x),\lambda_A(y))$ and $\lambda_2(x,y,z)=(\pi(x),\lambda_B(z)) $. It is easily seen that
$\Gamma(R/p\times A,R/p\times B)=\big\{((0)\times A,\;(0)\times B)\big\}$. By Theorem~\ref{z-main1}, the dimension of the tensor product of the $R\times R\times R$-algebras $((R/p\times A),\lambda_1)$ and $((R/p\times B),\lambda_2)$ is equal to 0. On the other hand, let $\lambda_2^{\prime}:R\times R\times R\rightarrow R/p\times B$ be a ring homomorphism defined by $\lambda_2^{\prime}(x,y,z)=(\pi(x),\lambda_B(y))$. By Theorem~\ref{z-main1}, the dimension of the tensor product of $((R/p\times A),\lambda_1)$ and $((R/p\times B),\lambda_2 ^{\prime})$
is equal to $\dim(A\otimes B)\neq 0$.
\end{example}

The next corollary handles the special case of domains.

%%%%%%%%%%%%%%%%%%%%%%%%%%%%%%%%%%%%%%%%%%%%%%%%%%%%%%%%%%%%%%%%%
%%%%%%%%%%%%%%%%%%%%%%%%%%%%%%%%%%%%%%%%%%%%%%%%%%%%%%%%%%%%%%%%%
\begin{corollary}
Let $(A,\lambda_A)$ be an AF-domain and $(B,\lambda_B)$ an algebra with $\dim(A\otimes B)\neq 0$. Set $t:=\td(A)$, $d:=\dim(A)$, and $p_{A}:=\Ker(\lambda_A)$. Then
$$\dim(A\otimes B)=\max\big\{\htt(\overline{Q}[t])+\min\big(t, d+\td(B/Q)\big)\mid \overline{Q}=\frac{Q}{p_{A}B}\in\Spec\big(\frac{B}{p_{A}B}\big)\big\}.$$
If, in addition, $B$ is a domain, then
$$\dim(A\otimes B)=\max\big\{\htt(Q[t])+\min\big(t, d+\td(B/Q)\big)\mid Q\in\Spec(B)\big\}.$$
\end{corollary}

The next main result extends Theorem~\ref{af-main2} to the zero-dimensional case, by establishing necessary and sufficient conditions for a tensor product of AF-rings to satisfy Wadsworth's formula on AF-domains over a field \cite[Theorem 3.8]{W}.

%%%%%%%%%%%%%%%%%%%%%%%%%%%%%%%%%%%%%%%%%%%%%%%%%%%%%%%%%%%%%%%%%
%%%%%%%%%%%%%%%%%%%%%%%%%%%%%%%%%%%%%%%%%%%%%%%%%%%%%%%%%%%%%%%%%
\begin{theorem}\label{z-main2}
Let  $(A_1,\lambda_1),\dots,(A_n,\lambda_n)$ be tensorially compatible AF-rings. Then, the following assertions are equivalent:
\begin{enumerate}[label= (\rm\roman*)]
\item $\dim(A_{1}\otimes\dots \otimes A_{n}) = \sum_{1\leq i\leq n}\td(A_{i})\ -\ \max\big\{\td(A_{i}) - \dim(A_{i})\big\}_{1\leq i\leq n}$;
\item For each $i =1,\dots, n$  there exists $M_i\in\Max(A_i)$ with $\lambda_{1}^{-1}(M_{1})=\dots=\lambda_{n}^{-1}(M_{n})$ such that $\htt(M_{i_{o}})=\dim(A_{i_{o}})$ for some $i_{o}\in\{1,..., n\}$ and, for all $i\not=i_{o}$, $\td(A_{iM_i}) = \td(A_{i})$ \& $\td(A_i/M_i) \leq\td(A_{i_{o}}/M_{i_{o}})$.
\end{enumerate}
\end{theorem}

%%%%%%%%%%%%%%%%%%%%%%%%%%%%%%%%%%%%%%%%%%%%%%%%%%%%%%%%%%%%%%%%%
%%%%%%%%%%%%%%%%%%%%%%%%%%%%%%%%%%%%%%%%%%%%%%%%%%%%%%%%%%%%%%%%%
\begin{corollary}
Let  $(A_1,\lambda_1),\dots,(A_n,\lambda_n)$ be tensorially compatible AF-rings. If anyone of the following conditions holds:
\begin{enumerate}[label= (\rm\alph*)]
\item For each $i =1,\dots, n$  there exists $M_i\in\Max(A_i)$ with $\lambda_{1}^{-1}(M_{1})=\dots=\lambda_{n}^{-1}(M_{n})$ such that $\htt(M_{i})=\dim(A_{i})$  and $\td(A_{iM_i}) = \td(A_{i})$
\item If $M_1,\dots,M_n$ are maximal ideals, respectively, of $A_1,\dots,A_n$  with $\lambda_{1}^{-1}(M_{1})=\dots=\lambda_{n}^{-1}(M_{n})$, then $\td(A_{iM_i}) = \td(A_{i})$ for $i=1,\dots,n$
\item If $P_1,\dots,P_n$ are minimal prime ideals, respectively, of $A_1,\dots,A_n$  with $\lambda_{1}^{-1}(P_{1})=\dots=\lambda_{n}^{-1}(P_{n})$, then $\td(A_{i}/{P_i}) = \td(A_{i})$ for $i=1,\dots,n$
\item $A_1,\dots,A_n$ are equicodimensional,
\end{enumerate}
then  $$\dim(A_{1}\otimes\dots \otimes A_{n}) = \sum_{1\leq i\leq n}\td(A_{i})-\max\big\{\td(A_{i}) - \dim(A_{i})\big\}_{1\leq i\leq n}.$$
\end{corollary}

The special case of AF-domains reads as follows.

%%%%%%%%%%%%%%%%%%%%%%%%%%%%%%%%%%%%%%%%%%%%%%%%%%%%%%%%%%%%%%%%%
%%%%%%%%%%%%%%%%%%%%%%%%%%%%%%%%%%%%%%%%%%%%%%%%%%%%%%%%%%%%%%%%%
\begin{corollary}
Let  $(D_1,\lambda_1),\dots,(D_n,\lambda_n)$ be tensorially compatible AF-domains. Then  $$\dim(D_{1}\otimes\dots \otimes D_{n}) = \sum_{1\leq i\leq n}\td(D_{i})-\max\big\{\td(D_{i}) - \dim(D_{i})\big\}_{1\leq i\leq n}.$$
\end{corollary}

The special case of $A\otimes A$ is given below.

%%%%%%%%%%%%%%%%%%%%%%%%%%%%%%%%%%%%%%%%%%%%%%%%%%%%%%%%%%%%%%%%%
%%%%%%%%%%%%%%%%%%%%%%%%%%%%%%%%%%%%%%%%%%%%%%%%%%%%%%%%%%%%%%%%%
\begin{corollary}
Let  $(A,\lambda_A)$  be an AF-ring. Then, the following assertions are equivalent:
\begin{enumerate}[label=(\rm\roman*)]
\item $\dim(A\otimes A) =\dim(A) + \td(A)$;
\item $\exists\ M,N\in\Max(A)$ with $\lambda_{A}^{-1}(M)=\lambda_{A}^{-1}(N)$ such that $\htt(M)=\dim(A)$, $\td(A_{N})=\td(A)$,  and $\td(A/N) \leq \td(A/M)$.
\end{enumerate}
\end{corollary}

%%%%%%%%%%%%%%%%%%%%%%%%%%%%%%%%%%%%%%%%%%%%%%%%%%%%%%%%%%%%%%%%%%%%%%%%%%%%%%%%%%%%%%%%%%%%%%%%%%%%%
%%%%%%%%%%%%%%%%%%%%%%%%%%%%%%%%%%%%%%%%%%%%%%%%%%%%%%%%%%%%%%%%%%%%%%%%%%%%%%%%%%%%%%%%%%%%%%%%%%%%%
%%%%%%%%%%%%%%%%%%%%%%%%%%%%%%%%%%%%%%%%%%%%%%%%%%%%%%%%%%%%%%%%%%%%%%%%%%%%%%%%%%%%%%%%%%%%%%%%%%%%%
%%%%%%%%%%%%%%%%%%%%%%%%%%%%%%%%%%%%%%%%%%%%%%%%%%%%%%%%%%%%%%%%%%%%%%%%%%%%%%%%%%%%%%%%%%%%%%%%%%%%%
\subsection{Transfer of the (locally) Jaffard property}

Theorem~\ref{af-main3} states that if $A$ is an AF-ring over a field $k$ and $B$ is a locally Jaffard $k$-algebra, then $A\otimes B$ is locally Jaffard. The main result of this subsection extends this result to AF-rings over a zero-dimensional ring.

%%%%%%%%%%%%%%%%%%%%%%%%%%%%%%%%%%%%%%%%%%%%%%%%%%%%%%%%%%%%%%%%%
%%%%%%%%%%%%%%%%%%%%%%%%%%%%%%%%%%%%%%%%%%%%%%%%%%%%%%%%%%%%%%%%%
\begin{theorem}\label{z-main3}
Let $A$ be an AF-ring (over $R$) and $B$ a locally Jaffard $R$-algebra with $A\otimes B\neq 0$. Then, $A\otimes B$ is locally Jaffard.
\end{theorem}

The next result asserts that Girolami's inequality on the valuative dimension \cite[Proposition 3.1]{Giro} holds in the zero-dimensional case.

%%%%%%%%%%%%%%%%%%%%%%%%%%%%%%%%%%%%%%%%%%%%%%%%%%%%%%%%%%%%%%%%%
%%%%%%%%%%%%%%%%%%%%%%%%%%%%%%%%%%%%%%%%%%%%%%%%%%%%%%%%%%%%%%%%%
\begin{proposition}
Let $A$ and $B$ be tensorially compatible algebras. Then
$$\dim_v(A_1\otimes_kA_2)\le\min\big(\dim_v(A_1) +\td(A_2),\td(A_1) + \dim_v(A_2)\big).$$
\end{proposition}

The next result handles the case where one of the two algebras is an AF-ring.

%%%%%%%%%%%%%%%%%%%%%%%%%%%%%%%%%%%%%%%%%%%%%%%%%%%%%%%%%%%%%%%%%
%%%%%%%%%%%%%%%%%%%%%%%%%%%%%%%%%%%%%%%%%%%%%%%%%%%%%%%%%%%%%%%%%
\begin{corollary}
Let $A$ be an AF-ring and $B$ an algebra with $\dim_{v}(B)\geq 1$ and $A\otimes B\neq 0$. Then, for any $r\geq \dim_{v}(B)-1$, we have
$$\dim_{v}(A\otimes B)= \max\big\{\htt(Q[r])+\min\big(\td(A_{P}),\htt(P)+\td(B/Q)\big)\mid (P,Q)\in\Gamma(A,B)\big\}.$$
\end{corollary}

If $A$ is an AF-domain, we get the following two results.

%%%%%%%%%%%%%%%%%%%%%%%%%%%%%%%%%%%%%%%%%%%%%%%%%%%%%%%%%%%%%%%%%
%%%%%%%%%%%%%%%%%%%%%%%%%%%%%%%%%%%%%%%%%%%%%%%%%%%%%%%%%%%%%%%%%
\begin{corollary}
Let $(A,\lambda_{A})$ be an AF-domain and $(B,\lambda_{B})$ an algebra with $\dim_{v}(B)\geq 1$ and $A\otimes B\neq 0$. Then, for any $r\geq \dim_{v}(B)-1$, we have
$$\dim_{v}(A\otimes B)= \max\big\{\htt(Q[r])+\min\big(\td(A),\dim(A)+\td(B/Q)\big)\big\}$$
where $Q$ ranges over the prime ideals of $B$ such that $\lambda_{B}^{-1}(Q)=\Ker(\lambda_{A})$.
\end{corollary}

%%%%%%%%%%%%%%%%%%%%%%%%%%%%%%%%%%%%%%%%%%%%%%%%%%%%%%%%%%%%%%%%%
%%%%%%%%%%%%%%%%%%%%%%%%%%%%%%%%%%%%%%%%%%%%%%%%%%%%%%%%%%%%%%%%%
\begin{corollary}
Let $A$ be an AF-domain and $B$ an algebra with $A\otimes B\neq 0$. If $\dim_{v}(B)\leq \td(A)+1$, then $A\otimes B$ is a Jaffard ring.
\end{corollary}

We conclude this section with the following observation. Let $A_{red}$ denote the reduced ring associated to a ring $A$. Then, $\td(A:R)=\td(A_{red}:R_{red})$, for any $R$-algebra $A$. Further,
if $(A,\lambda_A)$ and $(B,\lambda_B)$ are $R$-algebras, then
$(A\otimes_{R}B)_{red}=( A_{red}\otimes_{R_{red}}B_{red})_{red}$ by
\cite[Corollary 4.5.12]{GD}. Thus, one may assume that $R$ is absolutely flat
and $(A,\lambda_A)$ and $(B,\lambda_B)$ are reduced $R$-algebras.

%\newpage
%%%%%%%%%%%%%%%%%%%%%%%%%%%%%%%%%%%%%%%%%%%%%%%%%%%%%%%%%%%%%%%%%%%%%%%%%%%%%%%%%%%%%%%%%%%%%%%%%%%%%%%%%%%%%%%%%%%%%%%%%%%%%%%%%%%%%%%%%%%%%%
%%%%%%%%%%%%%%%%%%%%%%%%%%%%%%%%%%%%%%%%%%%%%%%%%%%%%%%%%%%%%%%%%%%%%%%%%%%%%%%%%%%%%%%%%%%%%%%%%%%%%%%%%%%%%%%%%%%%%%%%%%%%%%%%%%%%%%%%%%%%%%
%%%%%%%%%%%%%%%%%%%%%%%%%%%%%%%%%%%%%%%%%%%%%%%%%%%%%%%%%%%%%%%%%%%%%%%%%%%%%%%%%%%%%%%%%%%%%%%%%%%%%%%%%%%%%%%%%%%%%%%%%%%%%%%%%%%%%%%%%%%%%%
%%%%%%%%%%%%%%%%%%%%%%%%%%%%%%%%%%%%%%%%%%%%%%%%%%%%%%%%%%%%%%%%%%%%%%%%%%%%%%%%%%%%%%%%%%%%%%%%%%%%%%%%%%%%%%%%%%%%%%%%%%%%%%%%%%%%%%%%%%%%%%
\section{Tensor products of pullbacks issued from AF-domains}\label{p}

\noindent   This section is devoted to \cite{BGK3}, which establishes formulas for the Krull and valuative dimensions of tensor products of pullbacks issued from AF-domains. To this purpose, we use previous investigations of the prime ideal structure of various pullbacks, as in \cite{ABDFK,Ar,ArGi,BaGi,BK,BR,Cahen88,Cahen}. Moreover, in \cite{GiroKa}, a dimension formula for the tensor product of two particular pullbacks was established and a conjecture for more general pullbacks was raised; in this section, this conjecture is resolved.

Throughout, $k$ will be a field and $\mathcal{C}$ will denote the class of (commutative) $k$-algebras with finite transcendence degree over $k$. Algebras (resp., tensor products), when not specifically marked, will be taken over (resp., relative to) $k$.

Let $T$ be a domain, $M$ a maximal ideal of $T$, $K$ its residue field, $\varphi:T\longrightarrow K$ the canonical surjection, and $D$ a proper subring of $K$. Let $R$ be the pullback issued from the following diagram of canonical homomorphisms:
\[\begin{array}{cccl}
                    &R:=\varphi^{-1}(D) & \longrightarrow                       & D\\
(\ \square\ )       &\downarrow         &                                       & \downarrow\\
                    &T                  & \stackrel{\varphi}\longrightarrow     & K=T/M
\end{array}\]
Recall, from  \cite{Fon}, that $M = (R:T)$ and $D \cong R/M$; and for $p\in\Spec(R)$, if $M \not\subset p$, then $\exists !\ q\in\Spec(T)$ such that $q\cap R = p$ and $T_q = R_p$. However, if $M \subseteq p$, then $\exists !\ q\in\Spec(D)$ such that $p = \varphi^{-1}(q)$ and $R_{p}$ is a pullback determined by the following diagram
\[\begin{array}{cccl}
                    &R_{p}      & \longrightarrow                       & D_{q}\\
                    &\downarrow &                                       & \downarrow\\
                    &T_{M}      & \longrightarrow     & K
\end{array}\]
with $\htt(p)$ = $\htt(M)$ + $\htt(q)$. Recall also, from \cite{ABDFK,BR,Cahen88}, that
$$\dim(R) = \max\big\{\dim(T), \dim(D) + \dim(T_M)\big\}$$ and
$$\dim_v(R) = \max\big\{\dim_v(T), \dim_v(D) + \dim_v(T_M) + \t.d.(K:D)\big\}.$$ As for the dimension of the polynomial ring, we have the following lower bound which turned to be useful for the current study $$\dim(R[n])\geq \dim(D[n]) + \dim(T_M[n]) + \min\big(n, \t.d.(K:D)\big)$$
where the equality holds if $T$ is locally Jaffard with $\htt(M) = \dim(T)$.

%%%%%%%%%%%%%%%%%%%%%%%%%%%%%%%%%%%%%%%%%%%%%%%%%%%%%%%%%%%%%%%%%%%%%%%%%%%%%%%%%%%%%%%%%%%%%%%%%%%%%
%%%%%%%%%%%%%%%%%%%%%%%%%%%%%%%%%%%%%%%%%%%%%%%%%%%%%%%%%%%%%%%%%%%%%%%%%%%%%%%%%%%%%%%%%%%%%%%%%%%%%
%%%%%%%%%%%%%%%%%%%%%%%%%%%%%%%%%%%%%%%%%%%%%%%%%%%%%%%%%%%%%%%%%%%%%%%%%%%%%%%%%%%%%%%%%%%%%%%%%%%%%
%%%%%%%%%%%%%%%%%%%%%%%%%%%%%%%%%%%%%%%%%%%%%%%%%%%%%%%%%%%%%%%%%%%%%%%%%%%%%%%%%%%%%%%%%%%%%%%%%%%%%
\subsection{Krull dimension}

\noindent Recall that a pullback $R$ of type $\square$ is an AF-domain  if and only if $T$ and $D$ are AF-domains and $\t.d.(K:D) = 0$ \cite{Giro}. A combination of this result with the main result of this subsection  allows one to compute dimensions of tensor products for a large class of  algebras (that are not necessarily AF-domains).

The main theorem of this section relies on the following preliminaries which are important on their own. The next two lemmas deal with extensions of prime ideals of $R$ to polynomial rings over pullbacks.

%%%%%%%%%%%%%%%%%%%%%%%%%%%%%%%%%%%%%%%%%%%%%%%%%%%%%%%%%%%%%%%%%
%%%%%%%%%%%%%%%%%%%%%%%%%%%%%%%%%%%%%%%%%%%%%%%%%%%%%%%%%%%%%%%%%
\begin{lemma}
Let $R$ be a pullback of type $\square$ and $n$ a positive integer. For any $p\in\Spec(R)$ with $M\subseteq p$, we have
$$\htt(p[n])=\htt(p[n]/M[n])+\htt(M[n]).$$
\end{lemma}

%%%%%%%%%%%%%%%%%%%%%%%%%%%%%%%%%%%%%%%%%%%%%%%%%%%%%%%%%%%%%%%%%
%%%%%%%%%%%%%%%%%%%%%%%%%%%%%%%%%%%%%%%%%%%%%%%%%%%%%%%%%%%%%%%%%
\begin{lemma}
Let $R$ be a pullback of type $\square$ such that  $T_{M}$ and $D$ are locally Jaffard and let $n$ be a positive integer. For any $p\in\Spec(R)$ with $M\subseteq p$, we have
$$\htt(p[n])=\htt(p)+\min\big(n, \t.d.(K:D)\big).$$
\end{lemma}

The next two lemmas deal with the extensions of prime ideals to the tensor products.

%%%%%%%%%%%%%%%%%%%%%%%%%%%%%%%%%%%%%%%%%%%%%%%%%%%%%%%%%%%%%%%%%
%%%%%%%%%%%%%%%%%%%%%%%%%%%%%%%%%%%%%%%%%%%%%%%%%%%%%%%%%%%%%%%%%
\begin{lemma}
Let $A, B\in\mathcal{C}$ such that  $B$ is a domain. For any $p\in\Spec(A)$, we have
$$\htt(p\otimes B) = \htt(p[\td(B)]).$$
\end{lemma}

%%%%%%%%%%%%%%%%%%%%%%%%%%%%%%%%%%%%%%%%%%%%%%%%%%%%%%%%%%%%%%%%%
%%%%%%%%%%%%%%%%%%%%%%%%%%%%%%%%%%%%%%%%%%%%%%%%%%%%%%%%%%%%%%%%%
\begin{lemma}
Let $A, B\in\mathcal{C}$ such that  $B$ is an AF-domain. For any $P\in\Spec(A\otimes B)$ with $p:=P\cap A$, we have
$$\htt(P)=\htt(p\otimes B) + \htt\big(\frac{P}{p\otimes B}\big).$$
\end{lemma}

Let us fix notation for the rest of this section. Let $R_{1}$ and $R_{2}$ be two pullbacks of type $\square$ issued from the $k$-algebras ($\in\mathcal{C}$), respectively, $(T_{1},D_{1},K_{1}=T_{1}/M_{1})$ and $(T_{2},D_{2},K_{2}=T_{2}/M_{2})$. For $i=1,2$, set $d_{i}:=\dim(T_{i})$, $d'_{i}:=\dim(D_{i})$, $t_{i}:=\td(T_{i})$, $r_{i}:=\td(K_{i})$, and $s_{i}:=\td(D_{i})$. Also, we will use the functions $\delta(p,q)$ and $D(s,d,A)$ as defined in Section~\ref{af}.

%%%%%%%%%%%%%%%%%%%%%%%%%%%%%%%%%%%%%%%%%%%%%%%%%%%%%%%%%%%%%%%%%
%%%%%%%%%%%%%%%%%%%%%%%%%%%%%%%%%%%%%%%%%%%%%%%%%%%%%%%%%%%%%%%%%
\begin{lemma}
Assume $T_{1}$ and $T_{2}$ are AF-domains. For any $(p_{1},p_{2})\in\Spec(R_{1})\times\Spec(R_{2})$ with $M_{1}\not\subset p_{1}$ and $M_{2}\not\subset p_{2}$, we have
$$\delta(p_{1},p_{2})=\min\big(\htt(p_{1})+t_{2},t_{1}+\htt(p_{2})\big)\leq \min\big(d_{1}+t_{2},t_{1}+d_{2}\big).$$
\end{lemma}

%%%%%%%%%%%%%%%%%%%%%%%%%%%%%%%%%%%%%%%%%%%%%%%%%%%%%%%%%%%%%%%%%
%%%%%%%%%%%%%%%%%%%%%%%%%%%%%%%%%%%%%%%%%%%%%%%%%%%%%%%%%%%%%%%%%
\begin{lemma}
Assume $T_{1}$ and $T_{2}$ are AF-domains. For any $P\in\Spec(R_{1}\otimes R_{2})$ with $M_{1}\subseteq p_{1}:=P\cap R_{1}$ and $M_{2}\not\subset p_{2}:=P\cap R_{2}$, we have
$$\htt(P)=\htt(M_{1}[t_{2}]) + \htt\big(\frac{P}{M_{1}\otimes R_{2}}\big).$$
\end{lemma}

Next, we state the main theorem of this subsection.

%%%%%%%%%%%%%%%%%%%%%%%%%%%%%%%%%%%%%%%%%%%%%%%%%%%%%%%%%%%%%%%%%
%%%%%%%%%%%%%%%%%%%%%%%%%%%%%%%%%%%%%%%%%%%%%%%%%%%%%%%%%%%%%%%%%
\begin{theorem}\label{p-main1}
Assume $T_{1}$, $T_{2}$, $D_{1}$, and $D_{2}$ are AF-domains such that $\htt(M_{1})=d_{1}$ and $\htt(M_{2})=d_{2}$. Then
$$\dim(R_{1}\otimes R_{2})=\max\big\{\htt(M_{1}[t_{2}])+D(s_{1},d'_{1},R_{2}),\htt(M_{2}[t_{1}])+D(s_{2},d'_{2},R_{1})\big\}.$$
\end{theorem}

It is an open problem to compute $\dim(R_{1}\otimes R_{2})$ if only $T_1$ (or $T_2$) is
assumed to be an AF-domain. However, if both are not AF-domains, then
the above formula does not hold in general \cite[Examples 4.3]{W}.

The formula stated in the above theorem matches Wadsworth's formula in the particular case where $R_1$ and $R_2$ are AF-domains. Indeed, for $i:=1,2$, if  $R_i$ is an AF-domain, then so are $T_i$ and $D_i$ and $r_{i}=s_{i}$. Moreover, by \cite{ABDFK}, $n_{i}:=\dim(R_i)=d_{i}+d'_{i}$. So, the above theorem yields

$$\begin{array}{rl}
\dim(R_{1}\otimes R_{2})    &=\max\big\{\htt(M_{1}[t_{2}]) + \dim(D_{1}\otimes R_{2}),\htt(M_{2}[t_{1}]) + \dim(R_{1}\otimes D_{2})\big\}\\
                            &=\max\big\{d_{1}+\min(n_{2}+s_{1},t_{2}+d'_{1}),
                            d_{2}+\min(n_{1}+s_{2},t_{1}+d'_{2})\big\}\\
                            &=\max\big\{\min(n_{2}+r_{1}+d_{1},t_{2}+d'_{1}+d_{1}),
                            \min(n_{1}+r_{2}+d_{2},t_{1}+d'_{2}+d_{2})\big\}\\
                            &=\min(t_{1}+n_{2},t_{2}+n_{1}),\ \text{as desired}.
\end{array}$$

%%%%%%%%%%%%%%%%%%%%%%%%%%%%%%%%%%%%%%%%%%%%%%%%%%%%%%%%%%%%%%%%%%%%%%%%%%%%%%%%%%%%%%%%%%%%%%%%%%%%%
%%%%%%%%%%%%%%%%%%%%%%%%%%%%%%%%%%%%%%%%%%%%%%%%%%%%%%%%%%%%%%%%%%%%%%%%%%%%%%%%%%%%%%%%%%%%%%%%%%%%%
%%%%%%%%%%%%%%%%%%%%%%%%%%%%%%%%%%%%%%%%%%%%%%%%%%%%%%%%%%%%%%%%%%%%%%%%%%%%%%%%%%%%%%%%%%%%%%%%%%%%%
%%%%%%%%%%%%%%%%%%%%%%%%%%%%%%%%%%%%%%%%%%%%%%%%%%%%%%%%%%%%%%%%%%%%%%%%%%%%%%%%%%%%%%%%%%%%%%%%%%%%%
\subsection{Valuative dimension}

\noindent Recall for convenience, that the valuative dimension is stable under adjunction of indeterminates; i.e., $\dim_v(R[n]) = \dim_v(R) + n$, for any ring $R$ and any positive integer $n$ \cite[Theorem 2]{J}. However, the problem of computing the valuative dimension of the tensor product of two algebras is still elusively open. In \cite{Giro}, Girolami established a very useful upper bound for such an invariant; more exactly, she proved that if $A_1$ and $A_2$ are algebras, then $$\dim_v(A_1\otimes_kA_2)\le\min\big(\dim_v(A_1) +\td(A_2),\td(A_1) + \dim_v(A_2)\big).$$

The goal of this subsection is to compute the valuative dimension for a large
class of tensor products of algebras arising as pullbacks issued from AF-domains (and where the pullbacks are not necessarily AF-domains). To this purpose, the next two preliminary results establish the transfer of the notion of AF-domain to a polynomial ring over an arbitrary domain and over a pullback, respectively.

%%%%%%%%%%%%%%%%%%%%%%%%%%%%%%%%%%%%%%%%%%%%%%%%%%%%%%%%%%%%%%%%%
%%%%%%%%%%%%%%%%%%%%%%%%%%%%%%%%%%%%%%%%%%%%%%%%%%%%%%%%%%%%%%%%%
\begin{lemma}
Let $A$ be a domain ($\in\mathcal{C}$) and let $n$ be a positive integer. Then, $A[n]$ is an AF-domain if and only if $\htt(p[n])+\td(A/p)=\td(A),\ \forall\ p\in\Spec(A)$.
\end{lemma}

%%%%%%%%%%%%%%%%%%%%%%%%%%%%%%%%%%%%%%%%%%%%%%%%%%%%%%%%%%%%%%%%%
%%%%%%%%%%%%%%%%%%%%%%%%%%%%%%%%%%%%%%%%%%%%%%%%%%%%%%%%%%%%%%%%%
\begin{lemma}
Let $R$ be a pullback of type $\square$ such that  $T$ and $D$ are AF-domains. Then, the polynomial ring $R[\td(K)-\td(D)]$ is an AF-domain.
\end{lemma}

Next, we present  the main result of this subsection. Similarly to the previous subsection, we consider two pullbacks $R_{1}$ and $R_{2}$ of type $\square$ issued, respectively, from $(T_{1},D_{1},K_{1}=T_{1}/M_{1})$ and $(T_{2},D_{2},K_{2}=T_{2}/M_{2})$.

%%%%%%%%%%%%%%%%%%%%%%%%%%%%%%%%%%%%%%%%%%%%%%%%%%%%%%%%%%%%%%%%%
%%%%%%%%%%%%%%%%%%%%%%%%%%%%%%%%%%%%%%%%%%%%%%%%%%%%%%%%%%%%%%%%%
\begin{theorem}\label{p-main2}
Assume $T_{1}$, $T_{2}$, $D_{1}$, and $D_{2}$ are AF-domains such that $\htt(M_{1})=\dim(T_{1})$ and $\htt(M_{2})=\dim(T_{2})$. Then
$$\dim_{v}(R_{1}\otimes R_{2})=\min\big(\dim_{v}(R_{1})+\td(R_{2}),\dim_{v}(R_{2})+\td(R_{1})\big).$$
\end{theorem}

%%%%%%%%%%%%%%%%%%%%%%%%%%%%%%%%%%%%%%%%%%%%%%%%%%%%%%%%%%%%%%%%%%%%%%%%%%%%%%%%%%%%%%%%%%%%%%%%%%%%%
%%%%%%%%%%%%%%%%%%%%%%%%%%%%%%%%%%%%%%%%%%%%%%%%%%%%%%%%%%%%%%%%%%%%%%%%%%%%%%%%%%%%%%%%%%%%%%%%%%%%%
%%%%%%%%%%%%%%%%%%%%%%%%%%%%%%%%%%%%%%%%%%%%%%%%%%%%%%%%%%%%%%%%%%%%%%%%%%%%%%%%%%%%%%%%%%%%%%%%%%%%%
%%%%%%%%%%%%%%%%%%%%%%%%%%%%%%%%%%%%%%%%%%%%%%%%%%%%%%%%%%%%%%%%%%%%%%%%%%%%%%%%%%%%%%%%%%%%%%%%%%%%%
\subsection{Some applications and examples}

\noindent This subsection presents some applications of Theorem~\ref{p-main1} and Theorem~\ref{p-main2}. The first result features mild assumptions, on the transcendence degrees, for a tensor product of pullbacks issued from AF-domains to inherit the Jaffard property. As above, we consider two pullbacks $R_{1}$ and $R_{2}$ of type $\square$ issued, respectively, from $(T_{1},D_{1},K_{1}=T_{1}/M_{1})$ and $(T_{2},D_{2},K_{2}=T_{2}/M_{2})$; and, for $i=1,2$, we set $t_{i}:=\td(T_{i})$, $r_{i}:=\td(K_{i})$, and $s_{i}:=\td(D_{i})$.

%%%%%%%%%%%%%%%%%%%%%%%%%%%%%%%%%%%%%%%%%%%%%%%%%%%%%%%%%%%%%%%%%
%%%%%%%%%%%%%%%%%%%%%%%%%%%%%%%%%%%%%%%%%%%%%%%%%%%%%%%%%%%%%%%%%
\begin{theorem}
Assume $T_{1}$, $T_{2}$, $D_{1}$, and $D_{2}$ are AF-domains such that $M_{1}$ is the unique maximal ideal of $T_{1}$ with $\htt(M_{1})=\dim(T_{1})$ and $M_{2}$ is the unique maximal ideal of $T_{2}$ with $\htt(M_{2})=\dim(T_{2})$.  Then, the following assertions are equivalent:
\begin{enumerate}[label=(\rm\roman*)]
\item $R_{1}\otimes R_{2}$ is a Jaffard ring;
\item Either ``$r_{1}-s_{1}\leq t_{2}$ and $r_{2}-s_{2}\leq s_{1}$" or ``$r_{2}-s_{2}\leq t_{1}$ and $r_{1}-s_{1}\leq s_{2}$."
\end{enumerate}
\end{theorem}

The next result states, under weak assumptions, a formula for the Krull dimension similar to the one of Theorem~\ref{p-main1}. One may regard this result as an analogue of \cite[Theorem 5.4]{BaGi} (also \cite[Proposition 2.7]{ABDFK} and \cite[Corollary 1]{Cahen88}) in the special case of tensor products of pullbacks issued from AF-domains.

%%%%%%%%%%%%%%%%%%%%%%%%%%%%%%%%%%%%%%%%%%%%%%%%%%%%%%%%%%%%%%%%%
%%%%%%%%%%%%%%%%%%%%%%%%%%%%%%%%%%%%%%%%%%%%%%%%%%%%%%%%%%%%%%%%%
\begin{theorem}\label{p-main3}
Assume $T_{1}$, $T_{2}$, $D_{1}$, and $D_{2}$ are AF-domains such that $\htt(M_{1})=\dim(T_{1})$ and $\htt(M_{2})=\dim(T_{2})$. Suppose that either $s_{1}\leq r_{2}-s_{2}$ or $s_{2}\leq r_{1}-s_{1}$. Then
$$\dim(R_{1}\otimes R_{2})=\max\big\{\htt(M_{1}[t_{2}]) + \dim(D_{1}\otimes R_{2}),\htt(M_{2}[t_{1}]) + \dim(R_{1}\otimes D_{2})\big\}.$$
\end{theorem}

The next result handles the spacial case when $R_{1}=R_{2}$.

%%%%%%%%%%%%%%%%%%%%%%%%%%%%%%%%%%%%%%%%%%%%%%%%%%%%%%%%%%%%%%%%%
%%%%%%%%%%%%%%%%%%%%%%%%%%%%%%%%%%%%%%%%%%%%%%%%%%%%%%%%%%%%%%%%%
\begin{corollary}\label{p-cor1}
Let $R$ be a pullback of type $\square$ such that  $T$ is an AF-domain with $\htt(M)=\dim(T)$ and $D$ is a Jaffard domain. Set $t:=\td(T)$. Then
$$\dim(R\otimes R)=\htt(M[t]) + \dim(D\otimes R)$$
If, in addition, $\t.d.(K:D)\leq\td(D)$, then
$$\dim(R\otimes R)=\dim_{v}(R\otimes R)=t+\dim_{v}(R).$$
\end{corollary}

We close with some illustrative examples. The first example illustrates the fact that, in Theorem~\ref{p-main1} and Corollary~\ref{p-cor1}, the assumption ``$\htt(M_{i})=\dim(T_{i})$ $(i = 1, 2)$" is not superfluous.

%%%%%%%%%%%%%%%%%%%%%%%%%%%%%%%%%%%%%%%%%%%%%%%%%%%%%%%%%%%%%%%%%
%%%%%%%%%%%%%%%%%%%%%%%%%%%%%%%%%%%%%%%%%%%%%%%%%%%%%%%%%%%%%%%%%
\begin{example}
Let  $K$  be an algebraic extension of $k$, $T := S^{-1}K[X,Y]$, where $S := K[X,Y]\setminus \big((X)\cup (X - 1,Y)\big)$, and $M := S^{-1}(X)$. Consider the pullback $R$ of type $\square$ issued from $\big(T, k(Y),T/M=K(Y)\big)$. Since  $S^{-1}K[X,Y]$  is an AF-domain and $k(Y)\subset K(Y)$ is algebraic, then $R$  is an AF-domain \cite{Giro}. So,  $\dim(R\otimes R) = \dim(R) + \td(R) = 2 + 2 = 4$ by \cite[Corollary 4.2]{W}. Now, $\htt(M[2]) = \htt(M) = 1$ and  $\dim(k(Y)\otimes R) = \min(2, 1 + 2) = 2$. It follows that $\htt(M[2]) + \dim(k(Y)\otimes R) = 3\neq \dim(R\otimes R)$.
\end{example}

Next, we show how one can use Theorem~\ref{p-main1} to compute the Krull dimension of the tensor product of two algebras for a large class of algebras  (which are not necessarily AF-domains).

%%%%%%%%%%%%%%%%%%%%%%%%%%%%%%%%%%%%%%%%%%%%%%%%%%%%%%%%%%%%%%%%%
%%%%%%%%%%%%%%%%%%%%%%%%%%%%%%%%%%%%%%%%%%%%%%%%%%%%%%%%%%%%%%%%%
\begin{example}
Consider two pullbacks $R_{1}$ and $R_{2}$ of type $\square$ issued, respectively, from $\big(k(X,Y)[Z]_{(Z)},k(X),k(X,Y)\big)$ and $\big(k(X)[Z]_{(Z)},k,k(X)\big)$. We have Clearly, $\dim(R_1) = \dim(R_2) = 1$ and $\dim_v(R_1) =\dim_v(R_2) =2$. So, $R_1$ and $R_2$ are not AF-domains. By Theorem~\ref{p-main1}, $\dim(R_1\otimes R_2) = 4$. Now, notice that Wadsworth's formula fails here since $\min(\dim(R_1) + \td(R_2), \dim(R_2) + \td(R_1)) =3$.
\end{example}

A combination of Theorem~\ref{p-main1} and Theorem~\ref{p-main3} allows one to compute the Krull dimension of the tensor product for more general algebras, as shown by the next example.

%%%%%%%%%%%%%%%%%%%%%%%%%%%%%%%%%%%%%%%%%%%%%%%%%%%%%%%%%%%%%%%%%
%%%%%%%%%%%%%%%%%%%%%%%%%%%%%%%%%%%%%%%%%%%%%%%%%%%%%%%%%%%%%%%%%
\begin{example}
Consider two pullbacks $R_{1}$ and $R_{2}$ of type $\square$ issued, respectively, from $\big(k(X)[Y]_{(Y)},k,k(X)\big)$ and $\big(k(X,Y,Z)[T]_{(T)},R_{1},k(X,Y,Z)\big)$. We have, $\dim(R_{1})=1$ and $\dim_{v}(R_{1})=2$. So, $R_{1}$ is not an AF-domain and, by Theorem~\ref{p-main1}, we obtain $\dim(R_1\otimes R_1)=3$. Moreover,  $\dim(R_{2})=2$ and $\dim_{v}(R_{2})=4$. The conditions of Theorem~\ref{p-main1} do not hold for the pullbacks $R_{1}$ and $R_{2}$. We may, however, appeal to Theorem~\ref{p-main3} to get
$$\begin{array}{rl}
\dim(R_{1}\otimes R_{2})    &=\max\big\{\htt(M_{1}[4]) + \dim(k\otimes R_{2}),\htt(M_{2}[2]) + \dim(R_{1}\otimes R_{1})\big\}\\
                            &=\max\big\{2+2,2+3\big\}\\
                            &=5
\end{array}$$
where $M_{1}:=Yk(X)[Y]_{(Y)}$ and $M_{2}:=Tk(X,Y,Z)[T]_{(T)}$.
\end{example}

Next, we show how one can use Corollary~\ref{p-cor1} to construct examples of non-AF-domains $R$ where the tensor product $R\otimes R$ is Jaffard.

%%%%%%%%%%%%%%%%%%%%%%%%%%%%%%%%%%%%%%%%%%%%%%%%%%%%%%%%%%%%%%%%%
%%%%%%%%%%%%%%%%%%%%%%%%%%%%%%%%%%%%%%%%%%%%%%%%%%%%%%%%%%%%%%%%%
\begin{example}
Let $R$ be the pullback issued from $\big(k(X,Y,Z)[T]_{(T)},k(X,Y),k(X,Y,Z)\big)$. Clearly $R$ is not an AF-domain since $\dim(R)\neq \dim_{v}(R)$. Moreover, note that $\t.d.(k(X,Y,Z):k(X,Y))\lneqq\td(k(X,Y))$. By Corollary~\ref{p-cor1}, $\dim(R\otimes R)=\dim_{v}(R\otimes R)=5$. That is, $R\otimes R$ is a Jaffard ring.
\end{example}

%\newpage
%%%%%%%%%%%%%%%%%%%%%%%%%%%%%%%%%%%%%%%%%%%%%%%%%%%%%%%%%%%%%%%%%%%%%%%%%%%%%%%%%%%%%%%%%%%%%%%%%%%%%%%%%%%%%%%%%%%%%%%%%%%%%%%%%%%%%%%%%%%%%%
%%%%%%%%%%%%%%%%%%%%%%%%%%%%%%%%%%%%%%%%%%%%%%%%%%%%%%%%%%%%%%%%%%%%%%%%%%%%%%%%%%%%%%%%%%%%%%%%%%%%%%%%%%%%%%%%%%%%%%%%%%%%%%%%%%%%%%%%%%%%%%
%%%%%%%%%%%%%%%%%%%%%%%%%%%%%%%%%%%%%%%%%%%%%%%%%%%%%%%%%%%%%%%%%%%%%%%%%%%%%%%%%%%%%%%%%%%%%%%%%%%%%%%%%%%%%%%%%%%%%%%%%%%%%%%%%%%%%%%%%%%%%%
%%%%%%%%%%%%%%%%%%%%%%%%%%%%%%%%%%%%%%%%%%%%%%%%%%%%%%%%%%%%%%%%%%%%%%%%%%%%%%%%%%%%%%%%%%%%%%%%%%%%%%%%%%%%%%%%%%%%%%%%%%%%%%%%%%%%%%%%%%%%%%

\end{document}